\documentclass[preprint,sort&compress]{elsarticle}

\usepackage{overpic}
\usepackage{amsmath,amstext,amsbsy,amssymb,mathdots}
\usepackage{tikz}
\usetikzlibrary{positioning}
\usetikzlibrary{shapes,arrows}
\usepackage{amsfonts}
\usepackage{relsize}
\usepackage{psfrag}
\usepackage{enumerate}
\usepackage{url}
\usepackage{caption}

\newdefinition{rmk}{Remark}

\journal{JCP}

\usepackage{alltt}
\usepackage{color}
\usepackage{fullpage}
\definecolor{string}{rgb}{0.7,0.0,0.0}
\definecolor{comment}{rgb}{0,0.4,0}
\definecolor{keyword}{rgb}{0,0,1.0}

\newtheorem{remark}{Remark}

\newcommand{\beq}{\begin{equation}}
\newcommand{\eeq}{\end{equation}}
\newcommand{\beqs}{\begin{equation*}}
\newcommand{\eeqs}{\end{equation*}}
\newcommand{\beqa}{\begin{eqnarray}}
\newcommand{\eeqa}{\end{eqnarray}}
\newcommand{\beqas}{\begin{eqnarray*}}
\newcommand{\eeqas}{\end{eqnarray*}}

\newcommand{\ep}{\varepsilon}
\newcommand{\epz}{\ep\rightarrow 0}
\newcommand{\vx}{\mathbf{x}}
\newcommand{\vy}{\mathbf{y}}
\newcommand{\vxc}{\mathbf{\hat{x}}}
\newcommand{\uf}{\underline{f}}
\newcommand{\ug}{\underline{g}}

\newcommand{\vf}{\mathbf{f}}
\newcommand{\vs}{\mathbf{s}}
\newcommand{\vw}{\mathbf{w}}
\newcommand{\vr}{\mathbf{r}}

\newcommand{\vlam}{\boldsymbol{\lambda}}

\newcommand{\pade}{Pad\'{e}}

\newcommand{\cond}{\text{cond}}
\newcommand{\ds}{\displaystyle}

\newcommand{\matlab}{\textsc{Matlab\;}}

\usepackage{amsopn}

\def\calL{{ \mathcal D}}

\def \Xh{ \hat{X} } 
 
\def \Yh{ \hat{Y} }
\def \vyc{ \hat{\vy} } 
\def \tw{ \dot{w} } 

\begin{document}

\begin{frontmatter}
\title{Stable computations with flat radial basis functions using vector-valued rational approximations} 

\author[addr1]{Grady B. Wright\corref{corresp}}
\address[addr1]{Department of Mathematics, Boise State University, Boise, ID 83725-1555, USA}
\ead{gradywright@boisestate.edu}
\cortext[corresp]{Corresponding Author}

\author[addr2]{Bengt Fornberg}
\address[addr2]{Department of Applied  Mathematics, University of Colorado, 526 UCB, Boulder, CO 80309, USA }
\ead{Bengt.Fornberg@colorado.edu}

\begin{abstract}
One commonly finds in applications of smooth radial basis functions (RBFs) that scaling the kernels so they are `flat' leads to smaller discretization errors.  However, the direct numerical approach for computing with flat RBFs (RBF-Direct) is severely ill-conditioned.  We present an algorithm for bypassing this ill-conditioning that is based on a new method for rational approximation (RA) of vector-valued analytic functions with the property that all components of the vector share the same singularities. This new algorithm (RBF-RA) is more accurate, robust, and easier to implement than the Contour-Pad\'e method, which is similarly based on vector-valued rational approximation.  In contrast to the stable RBF-QR and RBF-GA algorithms, which are based on finding a better conditioned base in the same RBF-space, the new algorithm can be used with any type of smooth radial kernel, and it is also applicable to a wider range of tasks (including calculating Hermite type implicit RBF-FD stencils).  We present a series of numerical experiments demonstrating the effectiveness of this new method for computing RBF interpolants in the flat regime.  We also demonstrate the flexibility of the method by using it to compute implicit RBF-FD formulas in the flat regime and then using these for solving Poisson's equation in a 3D spherical shell.
\end{abstract}



\end{frontmatter}

\section{Introduction}

Meshfree methods based on smooth radial basis functions\ (RBFs) are finding increasing use in scientific computing as they combine high order accuracy with enormous geometric flexibility in applications such as interpolation and for numerically solving PDEs. In these applications, one finds that the best accuracy is often achieved when their shape parameter $\ep$ is small, meaning that they are relatively flat~\cite{FoPi08,LarFor03}. 

The so called Uncertainty Principle, formulated in 1995~\cite{Sch95}, has contributed to a widespread misconception that flat radial kernels unavoidably lead to numerical ill-conditioning. This `principle' mistakenly assumes that RBF interpolants need to be computed by solving the standard RBF linear system (often denoted RBF-Direct). However, it has now been known for over a decade\cite{DriscollFornberg2002,FoWrLa,LarFor05,Sch05} that that the ill-conditioning issue is specific to this RBF-Direct approach, and that it can be avoided using alternative methods.  Three distinctly different numerical algorithms have been presented thus far in the literature for avoiding this ill-conditioning and thus open up the complete range of $\ep$ that can be considered.  These are the Contour-\pade\; (RBF-CP) method~\cite{FoWr}, the RBF-QR method~\cite{FLF,FoPi07,FaMC12}, and the RBF-GA method~\cite{FoLePo13}.  The present paper develops a new stable algorithm that is in the same category as RBF-CP.

For a fixed number of interpolation nodes and evaluations points, an RBF interpolant can be viewed as a vector-valued function of $\ep$~\cite{FoWr}.  The RBF-CP method exploits the analytic nature of this vector-valued function in the complex $\ep$-plane to obtain a vector-valued rational approximation that can be used as a proxy for computing stably in the $\epz$ limit.  One key property that is utilized in this method is that all the components of the vector-valued function share the same singularities (which are limited to poles).  The RBF-CP method obtains a vector-valued rational approximation with this property from contour integration in the complex $\ep$-plane and Pad\'e approximation.  However, this method is somewhat computationally costly and can be numerically sensitive to the determination of the poles in the rational approximations.  In this paper, we follow a similar approach of generating vector-valued rational approximants, but use a newly developed method for computing these.  The advantages of this new method, which we refer to as RBF-RA, over RBF-CP include:
\begin{itemize}
\item Significantly higher accuracy for the same computational cost.
\item Shorter, simpler code involving fewer parameters, and less use of complex floating point arithmetic.
\item More robust algorithm for computing the poles of the rational approximation.
\end{itemize}

As with the RBF-CP method, the new RBF-RA method is limited to a relatively low number of interpolation nodes (just under a hundred in 2-D, a few hundred in 3-D), but is otherwise more flexible than RBF-QR and RBF-GA in that it immediately applies to any type of smooth RBFs (see Table \ref{tbl:RBFExamples} for examples), to any dimension, and to more generalized interpolation techniques, such as appending polynomials to the basis, Hermite interpolation, and customized matrix-valued kernel interpolation.  Additionally, it can be immediately applied to computing RBF generated finite difference formulas (RBF-FD) and Hermite (or compact or implicit) RBF-FD formulas (termed RBF-HFD), which are based on standard and Hermite RBF interpolants, respectively~\cite{Wright200699}.  RBF-FD formulas have seen tremendous applications to solving various PDEs~\cite{FFBook,FoL11,Sanyasiraju2008,ShuDing2003,SPLM,Wright200699,FlyerLehto2012,SWFKJSC2014,FWF14,FF15} since being introduced around 2002~\cite{Tolstykh2003,WangLiu2002}.  It is for computing RBF-FD and  RBF-HFD formulas that we see the main benefits of the RBF-RA method, as these formulas are typically based on node sizes well within its limitations.  Additionally, in the case of RBF-HFD formulas, the RBF-QR and RBF-GA methods cannot be readily used.  

Another two areas where RBF-RA is applicable is in the RBF partition of unity (RBF-PU) method~\cite{Wendland02,Safdari-Vaighani2015,Shcherbakov2016185,CaDeRoPe:2015} and domain decomposition~\cite{LH,ZHL}, as these also involve relatively small node sets.   While the RBF-QR and RBF-GA methods are also applicable for these problems, they are limited to the Gaussian (GA) kernel, whereas the RBF-RA method is not.  In the flat limit, different kernels sometimes give results of different accuracies. It is therefore beneficial to have stable algorithms that work for all analytic RBFs.  Figure 8 in~\cite{FoWr} shows an example where the flat limits of multiquadric (MQ) and inverse quadratic (IQ) interpolants are about two orders of magnitude more accurate than for GA interpolants.


The remainder of the paper is organized as follows.  We review the issues with RBF interpolation using flat kernels in Section \ref{sec:flatissues}.  We then discuss the new vector-valued rational approximation method that forms the foundation for the RBF-RA method for stable computations in Section \ref{sec:vvra}.  Section \ref{sec:analytic_nature} describes the analytic properties of RBF interpolants in the complex $\ep$-plane and how the new rational approximation method is applicable to computing these interpolants and also to computing RBF-HFD weights.  We present several numerical studies in Section \ref{sec:numerics}.  The first of these focuses on interpolation and illustrates the accuracy and robustness of the RBF-RA method over the RBF-CP method and also compares these methods to results using multiprecision arithmetic.  The latter part of Section \ref{sec:numerics} focuses on the application of the RBF-RA method to generating RBF-HFD formulas for the Laplacian and contains results from applying these formulas to solving Poisson's equation in a 3D spherical shell.  We make some concluding remarks about the method in Section \ref{sec:conclusion}.  Finally, a brief \matlab code is given in  \ref{appdx:code} and some suggestions on one main free parameters of the algorithms is given in \ref{appdx:contours}.



\section{The nature of RBF ill-conditioning in the flat regime}\label{sec:flatissues}

\begin{table}[tbp] 
\centering
\begin{tabular}[c]{|lll|}
\hline
\multicolumn{1}{|l}{\textbf{Name}} & \multicolumn{1}{|l}{\textbf{Abbreviation}} & \multicolumn{1}{|l|}{\textbf{Definition}}\\
\hline
Gaussian & GA & $\phi_{\ep}(r)=e^{-(\ep r)^{2}}$\\
Inverse quadratic & IQ & $\phi_{\ep}(r)=1/(1+(\ep r)^{2})$\\
Inverse multiquadric & IMQ & $\phi_{\ep}(r)=1/\sqrt{1+(\ep r)^{2}}$\\
Multiquadric & MQ & $\phi_{\ep}(r)=\sqrt{1+(\ep r)^{2}}$\\
\hline
\end{tabular}
\caption{Examples of analytic radial kernels featuring a shape-parameter $\ep$ that the RBF-RA procedure is immediately applicable to.  The first three kernels are positive-definite and the last is conditionally negative definite.
\label{tbl:RBFExamples}}
\end{table}

For notational simplicity, we will first focus on RBF interpolants of the form
\begin{equation}
s(\vx,\ep)=\sum_{i=1}^{N}\lambda_{i}\phi_{\ep}(\|\vx-\vxc_{i}\|), \label{s(x)}
\end{equation}
where $\{\vxc_i\}_{i=1}^{N}\subset\mathbb{R}^d$ are the interpolation nodes (or centers), $\vx\in \mathbb{R}^d$ is some evaluation point, $\|\cdot\|$ denotes the two-norm, and $\phi_{\ep}$ is an analytic radial kernel that is positive (negative) definite or conditionally positive (negative) definite (see Table \ref{tbl:RBFExamples} for common examples).  However, the RBF-RA method applies equally well to many other cases.  For example, if additional polynomial terms are included in (\ref{s(x)})~\cite{FlyerPHS}, if (\ref{s(x)}) is used to find the weights in RBF-FD formulas~\cite{FFBook}, or if more generalized RBF interpolants, such as divergence-free and curl-free interpolants~\cite{fuselier08}, are desired.  

\begin{figure}[tbp]
\centering
\includegraphics[width=0.7\textwidth]{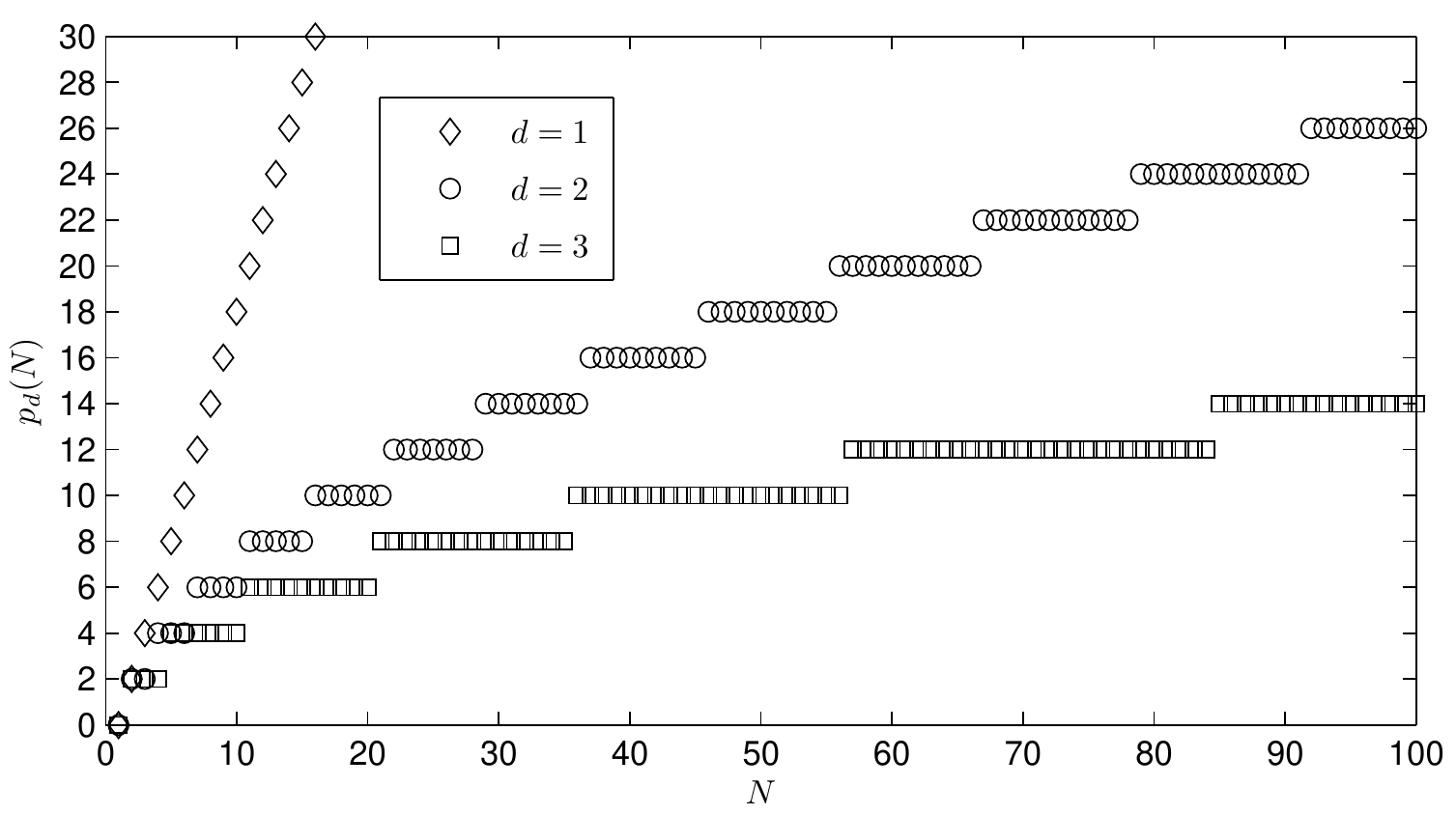}
\caption{The functions $p_{d}(N)$ in the condition number estimate $\cond(A(\ep))=O(\ep^{-p_{d}(N)})$ from~\cite{FZ07} for $d=1,2,3$.  The nodes are assumed to be scattered in $d$-D dimensions and the values of $p_{d}(N)$ are independent of the RBFs listed in Table \ref{tbl:RBFExamples}).}
\label{fig:sequences}
\end{figure}

The function $s(\vx ,\ep)$ interpolates the scattered data $\{\vxc_{i},g_{i}\}_{i=1}^N$ if the coefficients $\lambda_{i}$ are chosen as the solution of the linear system
\begin{equation}
A(\ep)\;\vlam(\ep)=\ug.
\label{eq:rbf_system}
\end{equation}
The matrix $A(\ep)$ has the entries
\begin{equation}
(A(\ep))_{i,j}=\phi_{\ep}(||\vxc_{i}-\vxc_{j}||)\;, \label{A-entries}
\end{equation}
and the column vectors $\vlam(\ep)$ and $\ug$ contain the $\lambda_{i}$ and the $g_{i}$ values, respectively.  We have explicitly indicated that the terms in \eqref{eq:rbf_system} depend on the choice of $\ep$ and we use underlines to indicate column vectors that do not depend on $\ep$ (otherwise they are bolded and explicitly marked). In the case of a fixed set of $N$ nodes scattered irregularly in $d$ dimensions, and when using any of the radial kernels listed in Table \ref{tbl:RBFExamples}, the condition number of the $A(\ep)$-matrix grows rapidly when the kernels are made increasingly flat (i.e $\epz$). As described first in~\cite{FZ07}, $\cond(A(\ep))=O(\ep^{-p_{d}(N)})$, where the functions $p_{d}(N)$ are illustrated in Figure \ref{fig:sequences}. The number of entries in the successive flat sections of the curves follow the pattern shown in Table \ref{tbl:sequences_table}.
\begin{table}[tbp] 
\centering
\begin{tabular}[c]{lrrrrrrr}
\hline
\hline
\textbf{Dimension} & \multicolumn{7}{l}{\textbf{Sequence}}\\
\hline
\hline
1-D & 1 & 1 & 1 & 1 & 1 & 1 & \ldots\\
2-D & 1 & 2 & 3 & 4 & 5 & 6 & \ldots\\
3-D & 1 & 3 & 6 & 10 & 15 & 21 & \ldots\\
\ldots &  &  &  &  &  &  & \\\hline\hline
\end{tabular}
\caption{The number of entries in the successive flat sections of the curves in Figure \ref{fig:sequences}.  If turned $45^{\circ}$ clockwise, this table coincides with Pascal's triangle. \label{tbl:sequences_table}}
\end{table}
Each row below the top one contains the partial sums of the previous row. As an example, with $N=100$ nodes (at the right edge of Figure \ref{fig:sequences}),
$\cond(A(\ep))$ becomes in 1-D $O(\ep^{-198}),$ in 2-D$\ O(\ep^{-26})$ and in 3-D $O(\ep^{-14})$.

As a result of the large condition numbers for small $\ep$, the $\lambda_{i}$-values obtained by (\ref{eq:rbf_system}) become extremely large in magnitude.  Since $s(\vx,\ep)$ depends in a perfectly well-conditioned way on the data $\{\vxc_{i},g_{i}\}_{i=1}^{N}$~\cite{DriscollFornberg2002,FoWrLa,Sch05}, a vast amount of numerical cancellation will then arise when the $O(1)$-sized quantity $s(\vx,\ep)$ is computed in (\ref{s(x)}). This immediate numerical implementation of (\ref{eq:rbf_system}) followed by (\ref{s(x)}) is known as the RBF-Direct approach. It consists of two successive ill-conditioned numerical steps for obtaining a well-behaved quantity. Like the RBF-CP, RBF-QR, and RBF-GA algorithms, the new RBF-RA method computes exactly the same quantity $s(\vx,\ep)$ by following steps that all remain numerically stable even in the $\epz$ limit. The condition numbers of the linear systems that arise within all these algorithms remain numerically reasonable even when $\ep\rightarrow0$. This often highly accurate parameter regime therefore becomes fully available for numerical work.

One strategy that has been attempted for coping with the ill conditioning of
RBF-Direct is to increase $\ep$ with $N$. If one keeps $\alpha=\ep/N^{1/d}$ constant in $d$-D, $\cond(A(\ep))$ stays roughly the same when $N$ increases. However, this approach causes \emph{stagnation (saturation) errors}~\cite{MazSchid,JPB1,FlyerPHS,Fasshauer:2007}, which severely damage the convergence properties of the RBF interpolant (and, for GA-type\ RBFs, causes convergence to fail altogether).  While convergence can be recovered by appending polynomial terms, it is reduced from spectral to algebraic.

\begin{figure}
\centering
\includegraphics[width=0.3\textwidth]{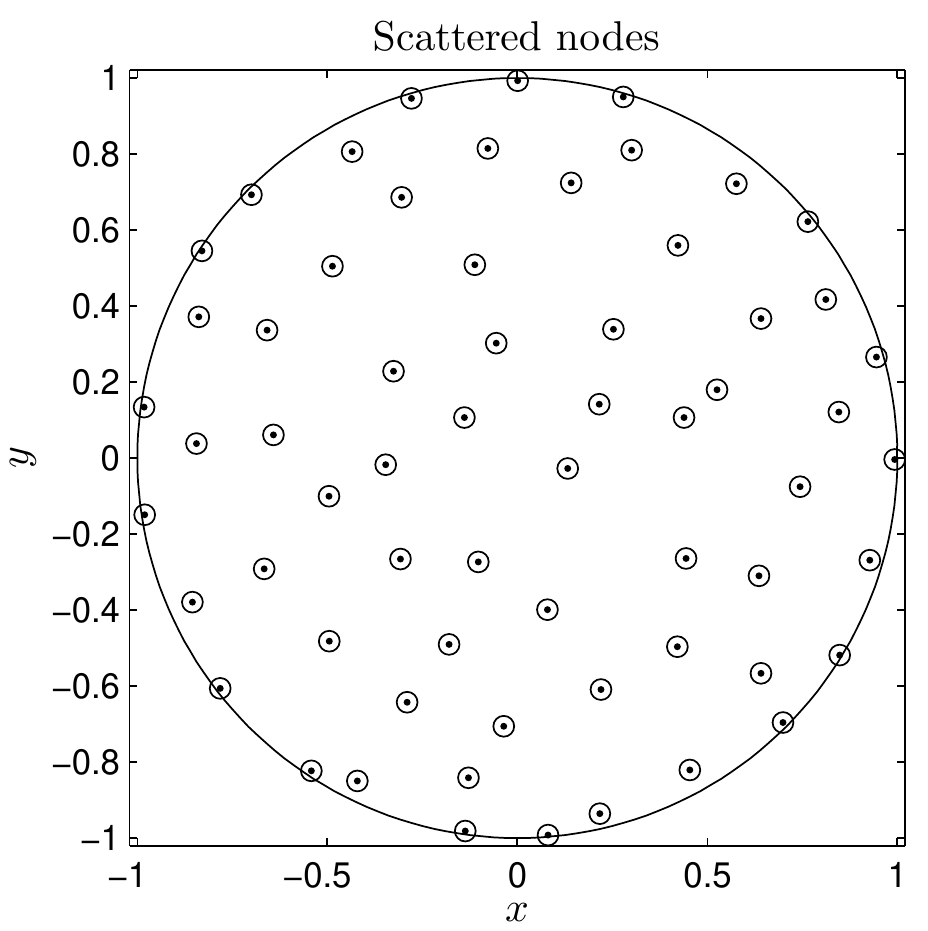}
\caption{Example of $N=62$ scattered nodes $\{\vxc_j\}_{j=1}^N$ over the unit circle.}
\label{fig:nodes}
\end{figure}

\section{Vector-valued rational approximation}\label{sec:vvra}
\begin{figure}
\centering
\includegraphics[width=0.33\textwidth]{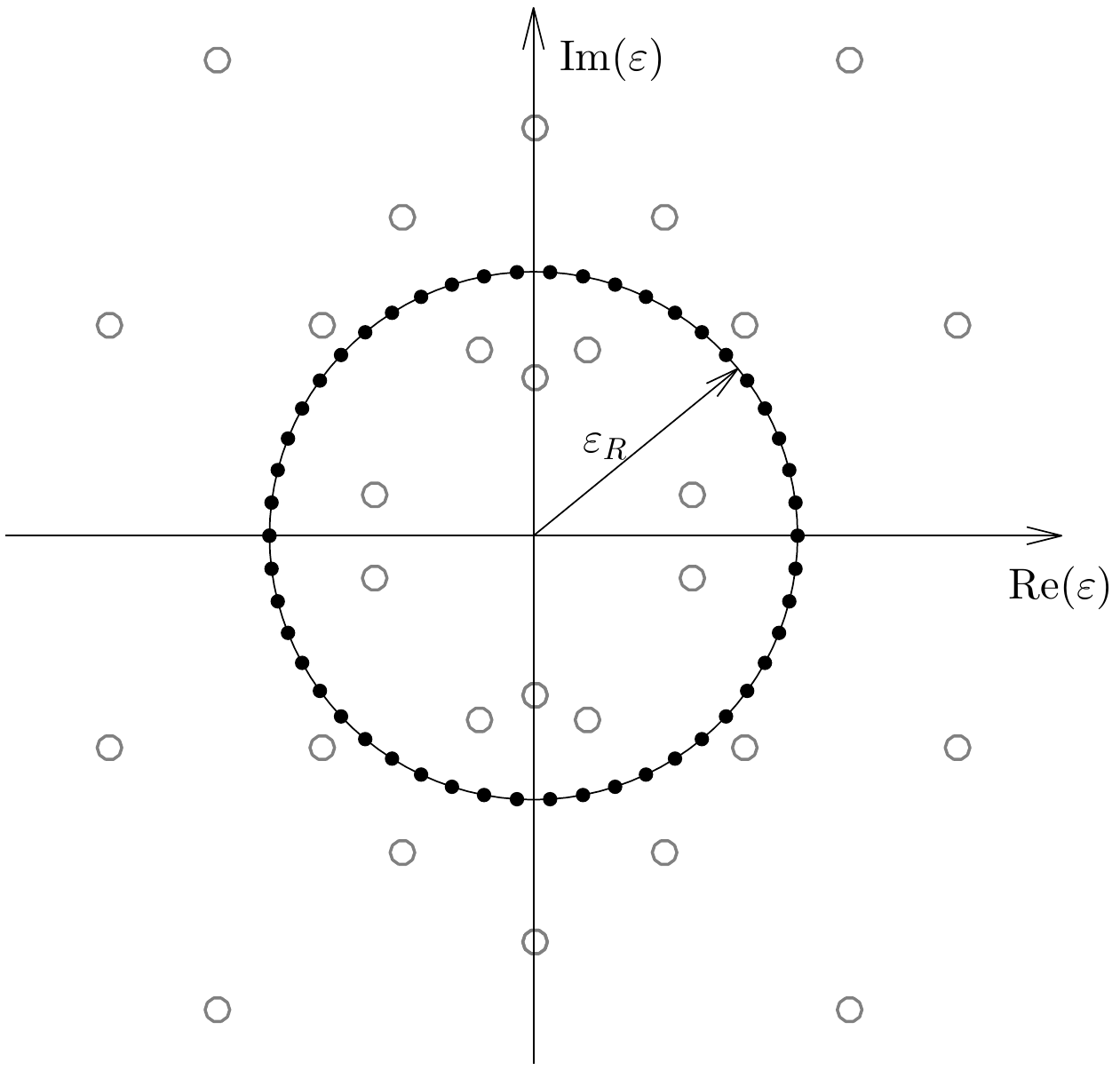}
\caption{Schematic of the analytic nature of the vector-valued functions $\vf(\ep)$ applicable to our vector-valued rational approximation method. The small open disks show the location of the poles common to all components of $\vf$ (taken here to be symmetric about the axes).  Solid black disks mark the locations where $\vf$ can be evaluated in a numerically stable manner.  The goal is to use these samples to construct a vector-valued rational approximation of the form \eqref{eq:rational} that can be used to accurately compute $\vf$ for $|\ep| < \ep_R$.\label{fig:schematicF}}
\end{figure}

The goal in this section is to keep the discussion of the vector-valued rational approximation method rather general since, as demonstrated in~\cite{Fornberg16}, it can be effective also for non-RBF-type applications.  Let $\vf(\ep):\mathbb{C}\rightarrow\mathbb{C}^M$ (with $M>1$) denote a vector-valued function with components $f_{j}(\ep)$, $j=1,\ldots,M$, that are analytic at all points within a region $\Omega$ around the origin of the complex $\ep$-plane except at possibly a finite number of isolated singular points (poles).  Furthermore, suppose that $\vf$ has the following properties:
\begin{enumerate}[(i)]
\item \vspace{5pt}All $M$-components of $\vf$ share the same singular points. 
\item Direct numerical evaluation of $\vf$ is only possible for $|\ep | \geq \ep_R > 0$,  where $|\ep | = \ep_R$ is within $\Omega$.
\item $\ep=0$ is at most a removable singular point of $\vf$.
\end{enumerate}
\vspace{5pt} We are interested in developing a vector-valued rational approximation to $\vf$ that exploits these properties and can be used to approximate $\vf$ for all $\ep < \ep_R$; see Figure \ref{fig:schematicF} for a representative scenario.  In Section \ref{sec:analytic_nature}, we return to RBFs and discuss how these approximations are especially applicable to the stable computation of RBF interpolation and  RBF-FD/HFD weights, both of which satisfy the properties (i)--(iii), for small $\ep$.  An additional property satisfied by these RBF applications is:
\begin{enumerate}[(i)]
\setcounter{enumi}{3}
\item \vspace{5pt}The function $\vf$ is even, i.e.\ $\vf(-\ep)=\vf(\ep)$.
\end{enumerate}  
\vspace{5pt} To simplify the discussion below, we will assume this property as well.   However, the present method is easily adaptable to the scenarios where this is not the case, and indeed, also to the case that property (iii) does not hold.

We seek to find a vector-valued rational approximation $\vr(\ep)$ to $\vf(\ep)$ with components taking the form
\begin{equation}
r_j(\ep) = \frac{a_{0,j}+a_{1,j}\ep^{2}+a_{2,j}\ep^{4}
+\ldots+a_{m,j}\ep^{2m}}{1+b_{1}\ep^{2}+b_{2}\ep
^{4}+\ldots+b_{n}\ep^{2n}}\,,\; j=1,\ldots,M. 
\label{eq:rational}
\end{equation}
It is important to note here that only the numerator coefficients depend on $j$ while the denominator coefficients are independent of $j$ to exploit property (i) above.  Additionally, we normalize the constant coefficient to one to match assumption (iii) and assume the numerators and denominator are even to match the assumption (iv).  To determine these coefficients, we first evaluate $\vf(\ep)$ around a circle of radius $\ep=\ep_R$ (see Figure \ref{fig:schematicF}), where $\vf$ can be numerically evaluated in a stable manner.  Due to assumption (iv), this evaluation only needs to be done at $K$ points $\ep_1,\ldots,\ep_K$ along the circle in the upper half-plane.  We then enforce that $r_j(\ep)$ agrees with $f_j(\ep)$ for all $K$ evaluation points to determine all the unknown coefficients of $\vr(\ep)$.  The enforcement of these conditions can, for each $j$, be written as the following coupled linear system of equations:
\begin{align}
\underbrace{
\begin{bmatrix}
 1 & \ep_1^2 & \cdots & \ep_1^{2m} \\
 1 & \ep_2^2 & \cdots & \ep_2^{2m} \\
 \vdots & \vdots & \ddots & \vdots \\ 
 1 & \ep_K^2 & \cdots & \ep_K^{2m}
\end{bmatrix}}_{\ds E}
\underbrace{
\begin{bmatrix}
a_{0,j} \\
\vdots \\
a_{m,j}
\end{bmatrix}}_{\ds \underline{a}_j}
+
\underbrace{\left(-\text{diag}(\underline{f}_j)
\begin{bmatrix}
\ep_1^2 & \cdots & \ep_1^{2n} \\
\ep_2^2 & \cdots & \ep_2^{2n} \\
\vdots & \ddots & \vdots  \\ 
\ep_K^2 & \cdots & \ep_K^{2n}
\end{bmatrix}
\right)}_{\ds F_j}
\underbrace{
\begin{bmatrix}
b_{1} \\
\vdots \\
b_{n}
\end{bmatrix}}_{\ds \underline{b}}
=
\underbrace{
\begin{bmatrix}
f(\ep_1) \\
f(\ep_2) \\
\vdots \\
f(\ep_K) \\
\end{bmatrix}}_{\ds \uf_j}.
\label{eq:coupled_system}
\end{align}
\begin{figure}[ptb]
\centering
\includegraphics[width=0.5\textwidth]{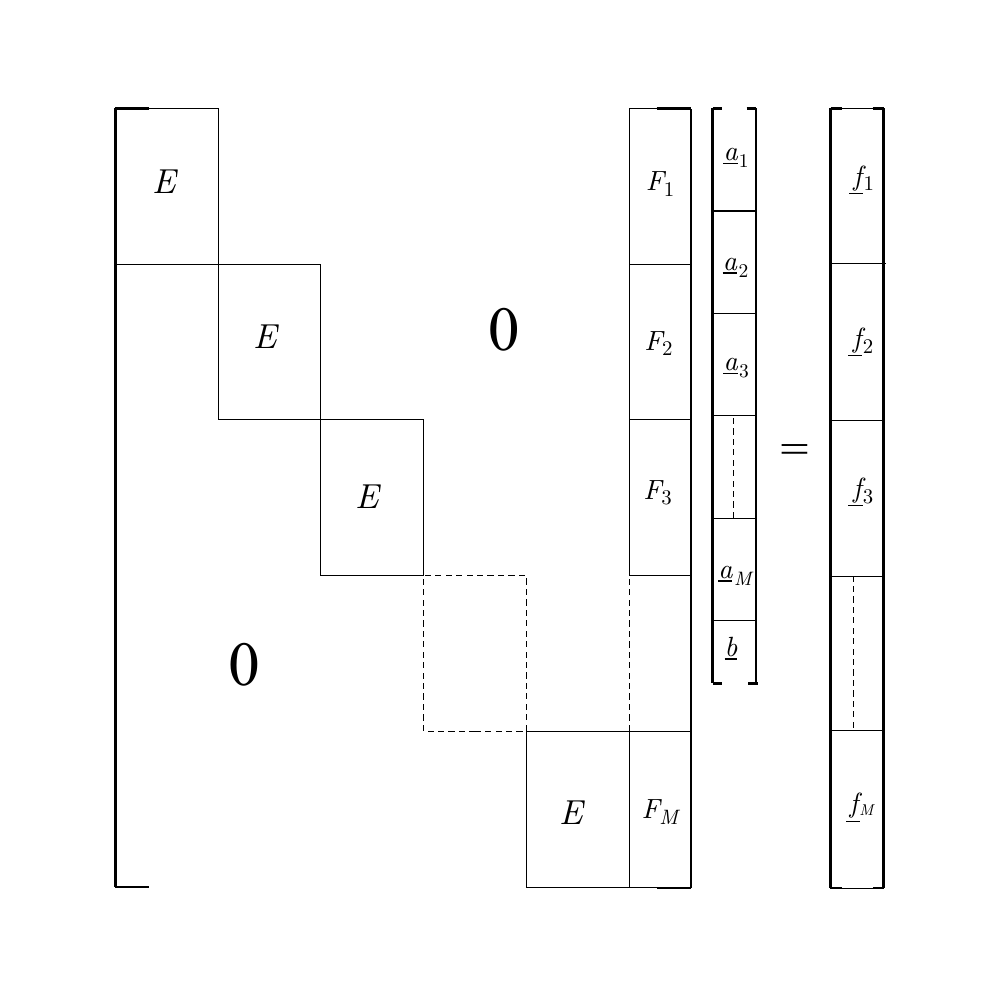}
\caption{Structure of overdetermined linear system \eqref{eq:coupled_system} for calculating the coefficients in the rational approximations (\ref{eq:rational}).\label{fig:system_1}}
\end{figure}
The structure of the complete system is displayed in Figure \ref{fig:system_1}.  In total there are $(m+1)M$ unknown coefficients corresponding to the different numerators of $\vr(\ep)$  and $n$ unknown coefficients corresponding to the one common denominator of $\vr(\ep)$.  Since each subsystem in \eqref{eq:coupled_system} is of size $K$-by-$K$, the complete system is of size $KM$-by-$((m+1)M+n)$.  We select parameters such that $m+1 + n/M <  K$ so that the system is overdetermined.

A schematic of the efficient algorithm we use to solve this coupled system is displayed in Figure \ref{fig:system_2}.  Below are some details on the five steps of this algorithm:
\begin{figure}[ptb]
\centering
\begin{tabular}{cc}
\parbox[t]{0.45\textwidth}{\centering \textbf{Step 1}: Normalize large rows among all block equations and compute a $QR$ factorization of $E$.} &
\parbox[t]{0.45\textwidth}{\centering \textbf{Step 2}: Multiply each block equation by $Q^{*}$.}\\
\includegraphics[width=0.45\textwidth]{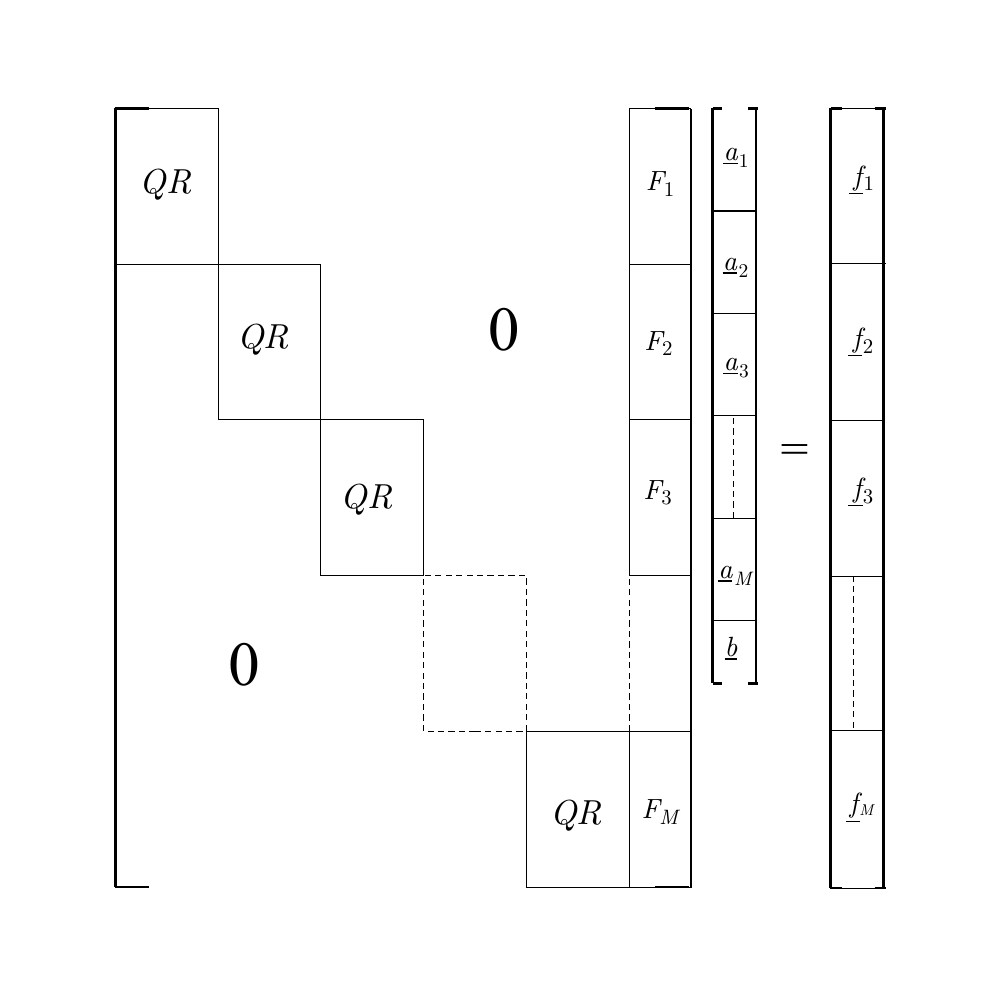} &
\includegraphics[width=0.455\textwidth]{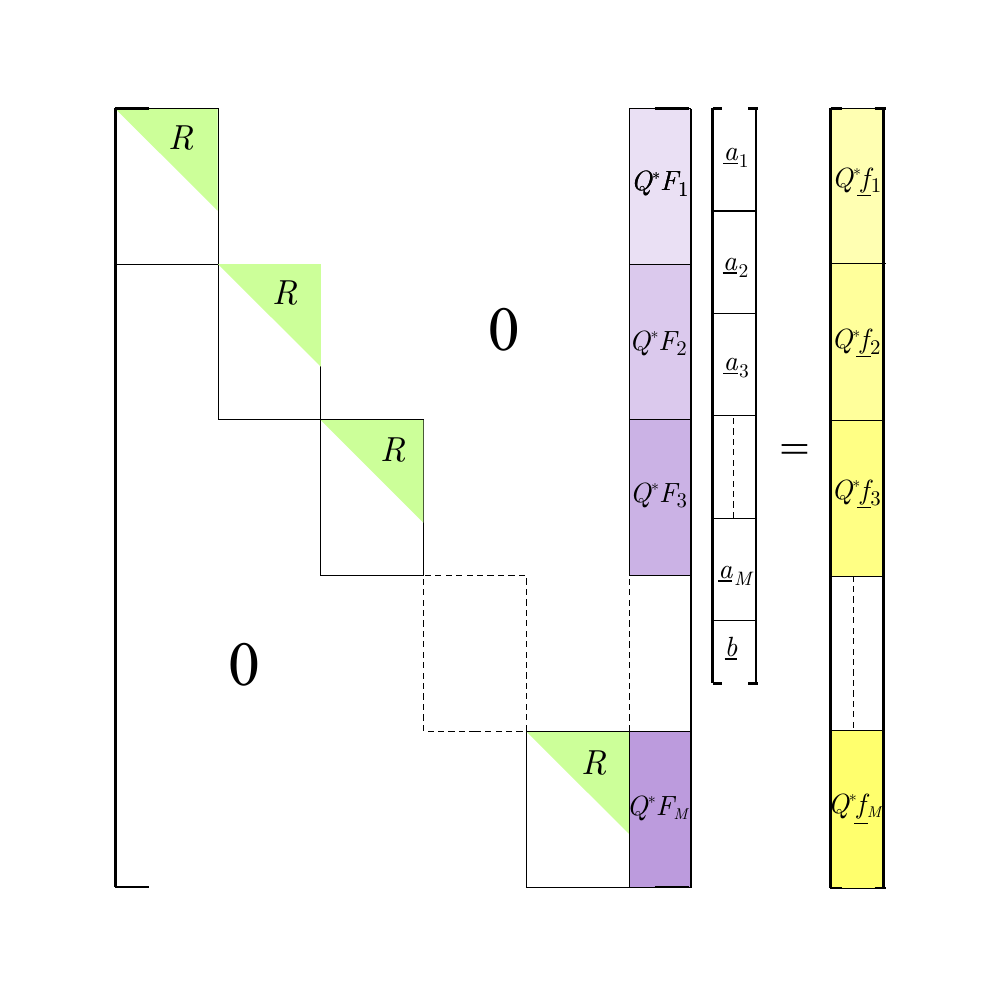} \\
\vspace{0.25\in} \\
\parbox[t]{0.45\textwidth}{\centering  \textbf{Step 3}: Re-order the rows so no zero-rows lie in the blocks on the diagonal.}  & \parbox[t]{0.45\textwidth}{\centering \textbf{Step 4}: Solve the decoupled, overdetermined system for $\underline{b}$ using least squares.}\\
\includegraphics[width=0.445\textwidth]{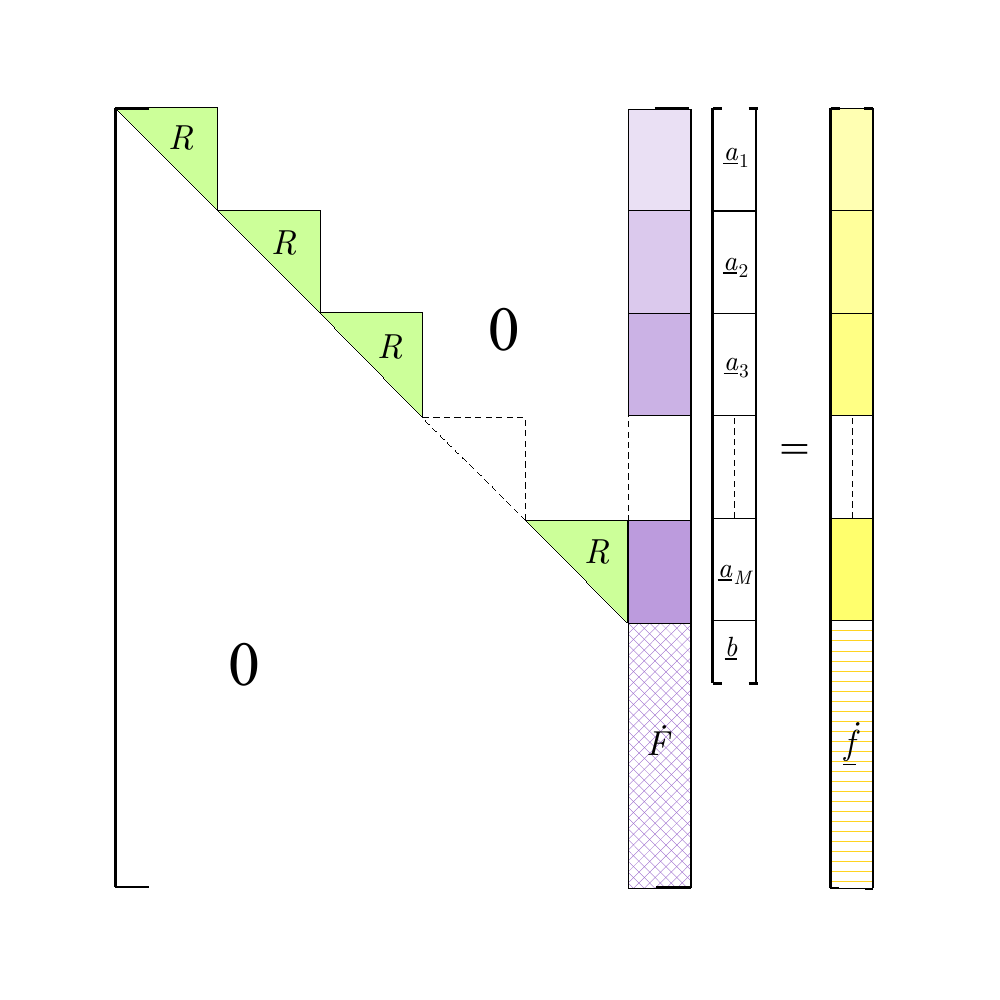} &
\parbox[t]{0.45\textwidth}{
\centering
\vspace{-2.9in}\includegraphics[width=0.2\textwidth]{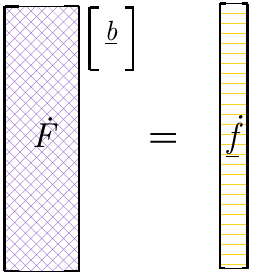} \\
\medskip
\textbf{Step 5}: Solve the remaining $M$ systems from Step 3 for $\underline{a}_j$ using $\underline{b}$ from Step 4.
\includegraphics[width=0.4\textwidth]{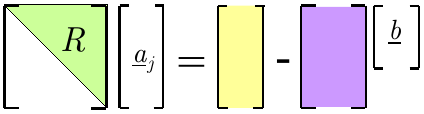}
}
\end{tabular}
\caption{Algorithm and schematic diagram for solving the system in Figure \ref{fig:system_1}.  Note that in the actual implementation, one does not actually need to construct the whole linear system and only one copy of $Q$ and $R$ is kept.
\label{fig:system_2}}
\end{figure}

\begin{enumerate} 
\item[Step 1:] \vspace{5pt} If an evaluation of $\vf$ at some $\ep_{k}$ happens to have been at or near a pole, then there may be an exceptionally large row in each $F_j$ matrix and corresponding row of $\uf_j$.  For square linear systems, multiplying a row by a constant will have no bearing on the system solution.  However, for least squares problems, it will influence the relative importance associated with the equation.  For each $\ep_k$, we thus normalize each of the associated rows (in $E$, $F_j$, and $\uf_j$) by dividing them by $\|\vf(\ep_k)\|_{\infty}$ to reduce their inflated importance in the system.  After this normalization, we compute a $QR$ factorization of the (modified) $E$ matrix.

\item[Step 2:] We left multiply the system by a block-diagonal matrix with blocks $Q^{*}$ on the diagonal.  This leaves the same upper triangular $R$ on the main diagonal blocks, with $Q^{*}F_j$ in last column of block matrices and $Q^{*}\uf_j$ in the right hand side blocks.

\item[Step 3:]  Starting from the top block equation, we re-order the rows so that the equations from the previous step corresponding to the zero rows of $R$ appear at the bottom of the linear system. This gives an almost upper triangular system, with a full matrix block of size $M(K-(m+1))$-by-$n$ in the last block column, which we denote by $\dot{F}$, and the corresponding rows of the right hand side, which we denote by $\dot{\uf}$.

\item[Step 4:] We compute the least squares solution to the $M(K-(m+1))$-by-$n$ overdetermined system $\dot{F}\underline{b} = \dot{\uf}$ for the coefficients of the common denominator of the vector-valued rational approximation.

\item[Step 5:] Using the coefficient vector $\underline{b}$ we finally solve the $M$ upper triangular block equations for the numerator coefficients $\underline{a}_j$, $j=1,\ldots,M$.  These systems are decoupled and consist of the same $(m+1)$-by-$(m+1)$ upper triangular block matrix $R$.
\end{enumerate}
A \matlab implementation of the entire vector-valued rational approximation method described above is given in \ref{appdx:code}, but with one modification that is applicable to RBFs.  We assume that $\vf(\ep) = \overline{\vf(\overline{\ep})}$, which additionally implies $\vf$ is real when $\ep$ is real.  This latter condition can be enforced on $\vr$ by requiring that the coefficients in the numerators and single denominator are all real.  Using this assumption and restriction on the coefficients of $\vr$, it then suffices to do the evaluations only along the circle of radius $\ep_{R}$ in the first quadrant and then split the resulting $K/2$ rows of the linear system \eqref{eq:coupled_system} into their real and imaginary parts to obtain $K$ rows that are now all real.  This allows the algorithm in Figure \ref{fig:system_2} to work entirely in real arithmetic.

The following are the dominant computational costs in computing the vector-valued rational approximation:  
\begin{enumerate} 
\item computing $f_j(\ep_k)$, $j=1,\dots,M$, and $k=1,\ldots,K,$ (or $k=1,\ldots,K/2,$ in the case that the coefficients are enforced to be real);
\item computing the $QR$ factorization of $E$, which requires $O(K(m+1)^2)$ operations;  
\item computing $Q^{*}F_j$, $j=1,\ldots,M$, which requires $O(MKmn)$ operations; 
\item solving the overdetermined system for $\underline{b}$, which requires $O(M(K-(m+1))n^2)$ operations; and
\item solving the $M$ upper triangular systems for $\underline{a}_j$, $j=1,\ldots,M$, which requires $O(M(m^2 + mn))$ operations.  
\end{enumerate}
We note that all of the steps of the algorithm are highly parallelizable, except for solving the overdetermined linear system for $\underline{b}$.  Additionally, we note that in cases where $M$ is so large that it dominates the cost, it may be possible to just use a subset of the evaluation points to determine $\underline{b}$, and then use these values to determine all of the numerator coefficients.  

In addition to choosing the radius of the contour $\ep_R$ for evaluating $\vf$, one also has to choose $m$, $n$, and $K$ in vector-valued rational approximation method.  We discuss these choices as they pertain specifically to RBF applications in the next section.

\begin{remark}
While we described the vector value rational approximation algorithm in terms of evaluation on a circular contour, this is not necessary.  It is possible to use more general evaluation contours, or even evaluations that are scattered within some band surrounding the origin where $\vf$ can be evaluated in a stable manner.
\end{remark}



\section{RBF-RA method}\label{sec:analytic_nature}
In this section, we describe how the vector-valued rational approximation method can be applied to both computing RBF interpolants and to computing RBF-FD/HFD weights.  We refer to the application of the vector-valued rational approximation method as the RBF-RA method.  We focus first on RBF interpolation since the insights gleaned  here can be applied to the RBF-FD/HFD setting (as well as many others).

\subsection{Stable computations of RBF Interpolants}\label{sec:interpolation}
Before the RBF-CP method was introduced, $\ep$ was in the literature always viewed as a real valued quantity. However, all the entries of the $A(\ep)$-matrix, as given by \eqref{A-entries}, are analytic functions of $\ep$ and, in the case of the GA kernel, they are entire functions of $\ep$. Thus, apart from any singularities associated with the kernel $\phi$ that is used to construct the interpolant, the entries of $A(\ep)^{-1}$ can have no other types of singularities than isolated poles.

\begin{figure}
\centering
\begin{tabular}{cc}
\includegraphics[width=0.55\textwidth]{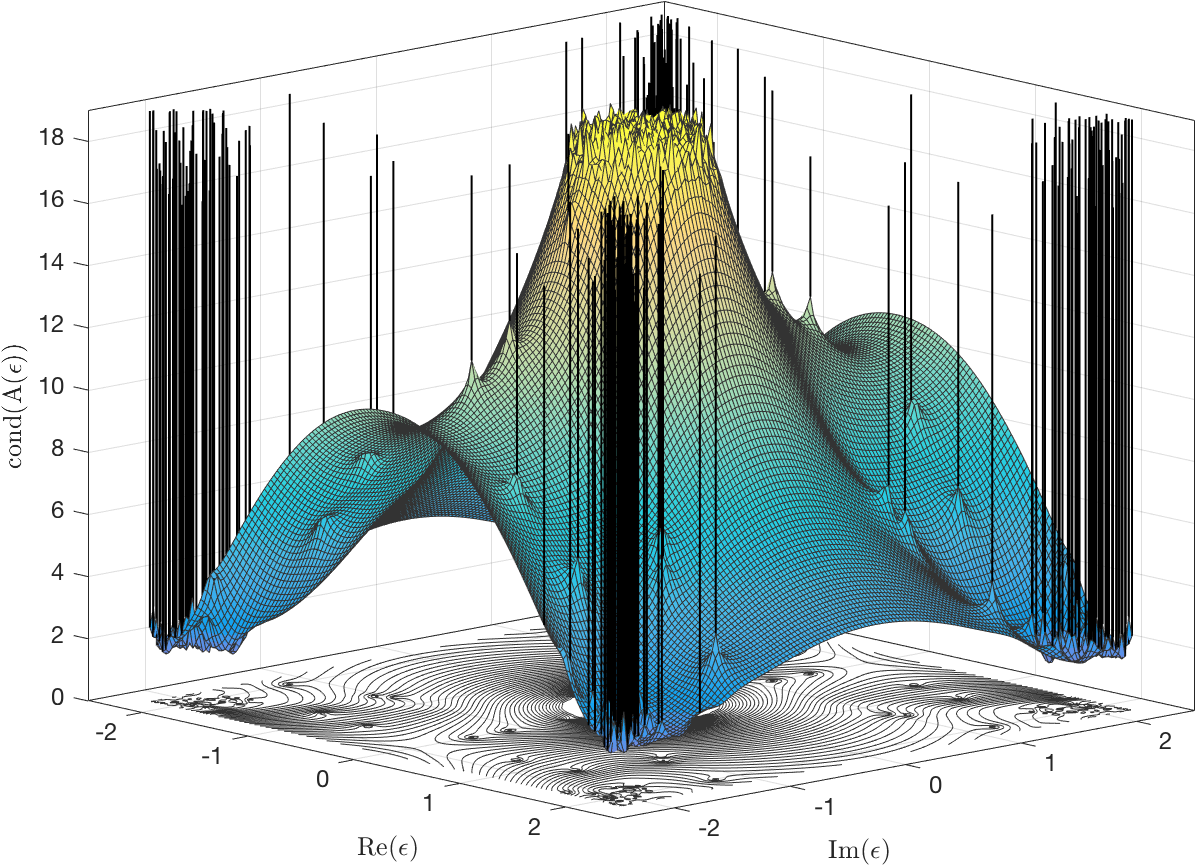} & 
\includegraphics[width=0.38\textwidth]{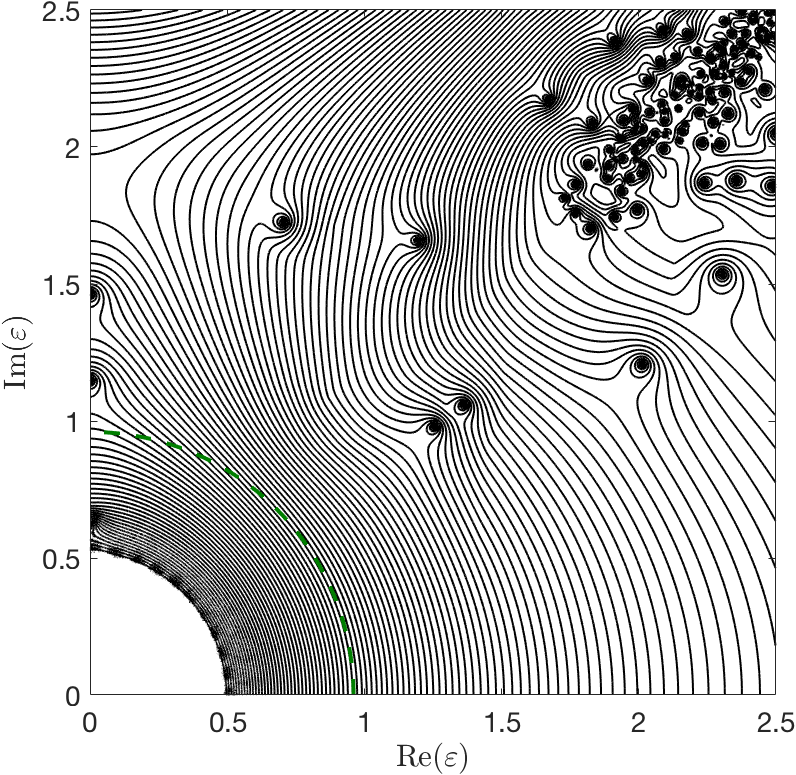} \\
(a) & (b)
\end{tabular}
\caption{ (a) Surface plot of $\log_{10}\cond(A(\ep))$ for the GA RBF and the $N=62$ nodes displayed in Figure \ref{fig:nodes} together with a contour plot with 51 equally spaced contours from 1 to 17.  (b) Contour plot of the surface in (a) over the first quadrant of the complex $\ep$-plane, now with 101 equally spaced contours from 1 to 17.  The dashed line marks a typical path along which the RBF-RA algorithm performs its RBF-Direct evaluations.}
\label{fig:cond_surface}
\end{figure}

To illustrate the behavior of $\cond(A(\ep))$ in the complex plane, consider the $N=62$ nodes $\{\vxc_{i}\}_{i=1}^N$ scattered over the unit circle in Figure \ref{fig:nodes} and the GA kernel.  For this example, $p_d(N)$ in Figure \ref{fig:sequences} is given by $p_{2}(62)=20$, so that $\cond(A(\ep))=||A(\ep)||_{2} ||A(\ep)^{-1}||_{2}=O(\ep^{-20})$. Since $||A(\ep)||_{2}$ stays $O(N)$ as $\epz$, it holds also that $||A(\ep)^{-1}||_{2}=O(\ep^{-20})$.  Figure \ref{fig:cond_surface} (a) shows a surface plot of $\log_{10}(\cond(A(\ep)))$ in the complex $\ep$-plane.
The first thing we note in this figure is that $\cond(A(\ep))$ is symmetric with respect to both the real and the imaginary axes since the interpolant (\ref{s(x)}) satisfies the symmetry relationship
\begin{align}
s(\vx,\ep) = s(\vx,-\ep) = \overline{s(\vx,\overline{\ep})} = \overline{s(\vx,-\overline{\ep})}.
\label{eq:Symmetry}
\end{align}
We also recognize the fast growth rate of $\cond(A(\ep))$ not only for $\epz$ along the real axis, but also when $\epz$ from any direction in the complex plane. Finally, we note sharp spikes, corresponding to the poles of $A(\ep)^{-1}$. To illustrate the behavior of $\cond(A(\ep))$ further, we display a contour plot of $\log_{10}(\cond(A(\ep)))$ in Figure \ref{fig:cond_surface} (b), but now only in the first quadrant, which is all that is necessary because of the symmetry conditions \eqref{eq:Symmetry}.
The number of poles is seen to increase very rapidly with increasing distance from the origin, but there are only a few present near the origin. The dashed line in Figure \ref{fig:cond_surface} (b) with radius $0.96$\ marks where, in this case $\cond(A(\ep))\approx10^{12}$, which is an acceptable level for computing $s(\vx,\ep)$ using RBF-Direct.

Now suppose that we are given $M>1$ points $\{\vx_j\}_{j=1}^M$ to evaluate the interpolant $s(\vx,\ep)$ at, with $\ep$ left as a variable.  We can write this as the following vector-valued function of $\ep$:
\begin{equation}
\underbrace{
\begin{bmatrix}
s(\vx_1,\ep) \\
s(\vx_2,\ep) \\
\vdots \\
s(\vx_M,\ep) 
\end{bmatrix}}_{\ds \vs(\ep)}
=
\underbrace{
\begin{bmatrix}
\phi_{\ep}(\|\vx_1-\vxc_{1}\|) & \cdots & \phi_{\ep}(\|\vx_1-\vxc_{N}\|) \\
\phi_{\ep}(\|\vx_2-\vxc_{1}\|) & \cdots & \phi_{\ep}(\|\vx_2-\vxc_{N}\|) \\
\vdots & \ddots & \vdots \\
\phi_{\ep}(\|\vx_M-\vxc_{1}\|) & \cdots & \phi_{\ep}(\|\vx_1-\vxc_{N}\|) \\
\end{bmatrix}}_{\ds \Phi(\ep)}
\begin{bmatrix}
& \vspace{7pt} & \\
& \phantom{A}A(\ep)^{-1}  & \\
& \vspace{7pt} & \\
\end{bmatrix}
\underbrace{
\begin{bmatrix}
g_{1}\vspace{5pt} \\
\vdots \\
g_{N}\vspace{5pt}
\end{bmatrix}}_{\ds \ug}.
\label{eq:rbf_mero}
\end{equation}
The entries of $\Phi(\ep)$ are analytic functions of $\ep$ within some disk of radius $\ep_R$ centered at the origin and $A(\ep)^{-1}$ can have at most poles in this disk.   Thus, $\vs(\ep)$ in \eqref{eq:rbf_mero} is analytic inside this disk apart from the poles of $A(\ep)^{-1}$.  Furthermore, since each entry of $\vs(\ep)$ is computed from $A(\ep)^{-1}$, they all share the same poles (note the similarity between Figure \ref{fig:cond_surface} (b) and Figure \ref{fig:schematicF}).  Finally, we know that typically, and always in the case that the GA kernel is used, that $\ep=0$ must be a removable singularity of $\vs(\ep)$~\cite{FoWrLa,Sch05}.  These observations together with the symmetry relationship \eqref{eq:Symmetry}, show that $\vs(\ep)$ satisfies all but possibly one of the properties for the vector-valued functions discussed in the previous section that the new vector-valued rational approximation is applicable to.  

The property that may not be satisfied of $\vs(\ep)$ is (ii).   The issue is that the region of ill-conditioning of $A(\ep)$, which tends to be disk-shaped, may spread out so far in the complex $\ep$-plane that it is not possible to find a radius $\ep_R$ for the evaluation contour where RBF-Direct can be used to solve \eqref{eq:rbf_mero} or that does not enclose the singular points associated with the kernels (branch points in the case of the MQ kernel and poles in the case of the IQ kernel).  The GA kernel has no singular points.  However, since it grows like $O(\exp(|\ep|^2))$ when $\pi/4 < \text{Arg}(\ep) < 3\pi/4$, this leads to a different type of ill-conditioning in $A(\ep)$ in these regions (see Figure \ref{fig:cond_surface}) and prevents a radius that is too large from being used.  For $N \lesssim 100$ in 2-D and $N \lesssim 300$ in 3-D, the region of ill conditioning surrounding $\ep=0$ is typically small enough (in double precision arithmetic) for a safe choice of $\ep_R$. As mentioned in the introduction, these limitations on $N$ are not problematic in many RBF applications.  We discuss some strategies for choosing $\ep_R$ in \ref{appdx:contours}.

After choosing the radius $\ep_R$ of the evaluation contour to use in the vector-valued rational approximation method for \eqref{eq:rbf_mero}, one must choose values for $m$, $n$, and $K$.   The optimal choice of these parameters depends on many factors, including the locations of the singular points from $A(\ep)^{-1}$, which are not known \emph{a priori}.  We have found that for a given $K$, choosing $m=K-1-n$ and $n = \lfloor K/4 \rfloor$ gives consistently good results across the problems we have tried.  With these choices, one can work out that the cost of the RBF-RA method is $O(KN(N^2+M) + MK^3)$.  As demonstrated by the numerical results in Section \ref{sec:numerics}, the approximations converge rapidly with $K$, so that not too large a value of $K$ is required to obtain an acceptable approximation.  With these parameters selected, we can construct a vector-valued rational approximation $\vr(\ep)$ for $\vs(\ep)$ in \eqref{eq:rbf_mero} that can be used as a proxy for computing the RBF interpolant stably for all values of $\ep$ inside the disk of radius $\ep_R$, including $\ep=0$.

\subsection{Stable computations of RBF-FD and RBF-HFD formulas}
RBF-FD and RBF-HFD formulas generalize standard, polynomial-based FD and compact or Hermite FD (HFD) formulas, respectively, to scattered nodes.  We discuss here only briefly how to generate these formulas and how they can be computed using the vector-valued rational approximation scheme from Section \ref{sec:vvra}.  More thorough discussions of the RBF-FD methods can be found in~\cite{FFBook}, and more details on RBF-HFD methods can be found in~\cite{Wright200699}.

Let $\calL$ be a differential operator (e.g.\ the Laplacian) and $\Xh = \{\vxc_i\}_{i=1}^N$ denote a set of (scattered) node locations in $\mathbb{R}^d$.  Suppose $g:\mathbb{R}^d\rightarrow\mathbb{R}$ is some (sufficiently smooth) function sampled on $\Xh$ and that we wish to approximate $\calL g(\vx)$ at $\vx=\vxc_1$ using a linear combination of samples of $g$ at the nodes in $\Xh$, i.e.\
\begin{align}
\calL g(\vx)\bigl |_{\vx=\vxc_1} \approx \sum_{i=1}^{N} w_i g(\vxc_i).
\label{eq:rbf_fd}
\end{align}
The RBF-FD method determines the weights $w_i$ in this approximation by requiring that \eqref{eq:rbf_fd} be exact whenever $g(\vx) = \phi_\ep(\|\vx-\vxc_i\|)$, $i=1,\ldots,N$, where $\phi_\ep$ is, for example, one of the radial kernels from Table \ref{tbl:RBFExamples}.\footnote{The RBF-FD (and HFD) method is general enough to allow for radial kernels that are not analytic functions of $\ep$, or that do not depend on a shape parameter at all~\cite{FlyerPHS}. We assume the kernels are analytic in the presentation as the RBF-RA algorithm is applicable only in this case.}  These conditions can be written as the following linear system of equations:
\begin{equation}
\begin{bmatrix}
& \vspace{9pt} & \\
& \phantom{A}A(\ep)\phantom{A}  & \\
& \vspace{9pt} & 
\end{bmatrix}
\underbrace{
\begin{bmatrix}
w_1(\ep) \\
w_2(\ep) \\
\vdots \\
w_N(\ep) 
\end{bmatrix}}_{\ds \vw(\ep)}
=
\underbrace{
\begin{bmatrix}
\calL_\vx\phi_{\ep}(\|\vx-\vxc_{1}\|)\bigl |_{\vx=\vxc_1} \\
\calL_\vx\phi_{\ep}(\|\vx-\vxc_{2}\|)\bigl |_{\vx=\vxc_1} \\
\vdots \\
\calL_\vx\phi_{\ep}(\|\vx-\vxc_{N}\|)\bigl |_{\vx=\vxc_1}
\end{bmatrix}}_{\ds \calL_\vx \boldsymbol{\phi}_\ep},
\label{eq:rbf_fd_mero1}
\end{equation}
where $A(\ep)$ is the standard RBF interpolation matrix with entries given in \eqref{A-entries} and $\calL_{\vx}$ means $\calL$ applied with respect to $\vx$ (this latter notation aids the discussion of the RBF-HFD method below).  We have here explicitly marked the dependence of the weights $w_i$ on $\ep$.  Using the same arguments as in the previous section, we see that the vector of weight $\vw(\ep)$ can be viewed as a vector-valued function of $\ep$ with similar analytic properties as an RBF interpolant based on $\phi_{\ep}$.  We can thus similarly apply the vector-valued rational approximation algorithm for computing $\vw(\ep)$ in a stable manner as $\epz$.

The additional constraint that \eqref{eq:rbf_fd} is exact when $g$ is a constant is often imposed in the RBF-FD method, which is equivalent to requiring $\sum_{i=1}^{N} w_i = \calL 1$.  This leads to an equality constrained quadratic programing problem that can be solved using Lagrange multipliers leading to the following slightly modified version of the system in~\eqref{eq:rbf_fd_mero1} for determining $\vw(\ep)$:
\begin{equation}
\begin{bmatrix}
A(\ep) & \mathbf{e} \\
\mathbf{e}^{T} & 0 
\end{bmatrix}
\begin{bmatrix}
\vw(\ep) \\
\lambda(\ep)
\end{bmatrix}
=
\begin{bmatrix}
\calL_{\vx} \boldsymbol{\phi}_\ep \\
\calL 1
\end{bmatrix},
\label{eq:rbf_fd_mero2}
\end{equation}
where $\mathbf{e}$ is the vector of size $N$ containing all ones, and $\lambda(\ep)$ is the Lagrange multiplier.  The vector-valued rational approximation algorithm is equally applicable to approximating $\vw(\ep)$ in this equation.

The HFD method is similar to the FD method, except that we seek to find an approximation of $\calL g(\vx)$ at $\vx=\vxc_1$ that involves a linear combination of $g$ at the nodes in $\Xh$ \textit{and} a linear combination of $\calL g$ at another set of nodes $\Yh = \{\vyc_j\}_{j=1}^L$, i.e.
\begin{align}
\calL g(\vx)\bigl |_{\vx=\vxc_1} \approx  \sum_{i=1}^{N} w_i g(\vxc_i) +  \sum_{j=1}^{L} \tw_j \calL g(\vx)\bigl |_{\vx=\vyc_j}.
\label{eq:rbf_hfd}
\end{align}
Here $\vxc_1\notin\Yh$, or else the approximation is trivial, also $\Yh$ is typically chosen to be some subset of $\Xh$ (not including $\vxc_1$), but this is not necessary.  The RBF-HFD method from~\cite{Wright200699} determines the weights by requiring that \eqref{eq:rbf_hfd} is exact for $g(\vx) = \phi_{\ep}(\|\vx-\vxc_i\|)$, $i=1,\ldots,N$, in addition to $g(\vx) = \calL_{\vy}  \phi_{\ep}(\|\vx-\vy\|)\bigl |_{\vy=\vyc_j}$, where $\calL_{\vy}$ means $\calL$ applied to the variable $\vy$.  This gives $N+L$ conditions for determining the weights in \eqref{eq:rbf_hfd}.  We can write these conditions as the following linear system for the vector of weights:
\begin{align}
\underbrace{
\begin{bmatrix}
A(\ep) & B(\ep) \\
B(\ep)^T & C(\ep)
\end{bmatrix}}_{\ds \tilde{A}(\ep)}
\underbrace{
\begin{bmatrix}
\vw(\ep) \\
\mathbf{\tw}(\ep)
\end{bmatrix}}_{\ds \tilde{\vw}(\ep)}
=
\begin{bmatrix}
\calL_{\vx} \boldsymbol{\phi}_{\ep} \\
\calL_{\vx}\calL_{\vy} \boldsymbol{\phi}_{\ep}
\label{eq:rbf_hfd_mero}
\end{bmatrix},
\end{align}
where $\tilde{\vw}(\ep) = \begin{bmatrix} w_1(\ep) & \cdots & w_N(\ep) & \tw_1(\ep) & \cdots & \tw_L(\ep)\end{bmatrix}^T$, $A(\ep)$ is the standard interpolation matrix for the nodes in $\Xh$ (see \eqref{A-entries}),
\begin{align*}
(B(\ep))_{ij} &= \calL_{\vy} \phi(\|\vx_i - \vy\|)\bigl |_{\vy = \vyc_j},\; i=1,\ldots,N,\; j=1,\ldots,L, \\
(C(\ep))_{ij} &= \calL_{\vx} \left(\calL_{\vy} \phi(\|\vx - \vy\|)\bigl |_{\vy = \vyc_j}\right)\bigl |_{\vx=\vxc_i},\; i=1,\ldots,L,\; j=1,\ldots,L,\text{and}\\
(\calL_{\vx}\calL_{\vy} \boldsymbol{\phi}_{\ep})_{i} &= \calL_{\vx} \left(\calL_{\vy} \phi(\|\vx - \vy\|)\bigl |_{\vy = \vyc_i}\right)\bigl |_{\vx=\vxc_1},\; i=1,\ldots,L.
\end{align*}
The system \eqref{eq:rbf_hfd_mero} is symmetric and, under very mild restrictions on the operator $\calL$, is non-singular for all the kernels $\phi_{\ep}$ in Table \ref{tbl:RBFExamples}; see~\cite{Wu1992} for more details.

Similar to the RBF-FD method, the vector of weights $\tilde{\vw}(\ep)$ in \eqref{eq:rbf_hfd_mero} is an analytic vector-valued function of $\ep$ with each entry sharing the same singular points, due now to $\tilde{A}(\ep)^{-1}$, in a region surrounding the origin.  We can thus similarly apply the vector-valued rational approximation algorithm for computing $\tilde{\vw}(\ep)$ in a stable manner as $\epz$. 

We note that in the RBF-HFD method it is also common to enforce that \eqref{eq:rbf_hfd} is exact for constants.  This constraint can be imposed in a likewise manner as the RBF-FD case in \eqref{eq:rbf_fd_mero2}.  We thus skip the details and refer the reader to~\cite{Wright200699} for the explicit construction.

\begin{remark}
In using the vector-valued rational approximation algorithm for RBF applications, we have found that it is beneficial to spatially re-scale the node sets by the radius of the evaluation contour in the complex $\ep$-plane, which allows all evaluations and subsequent computations of the algorithm to be done on the unit circle.  This is the approach we follow in the examples given in \ref{appdx:code}.
\end{remark}

\section{Numerical results}\label{sec:numerics}
In this section, we first present results of the RBF-RA algorithm as applied to RBF interpolation and follow this up with some results on computing RBF-HFD formulas and their application to a `large-scale' problem.  The primary studies of the accuracy and robustness of the algorithm in comparison to the RBF-CP method and multi-precision arithmetic are given for RBF interpolation as the observations made here also carry over to the RBF-HFD setting.

\subsection{RBF Interpolation results}\label{sec:interp_results}

\begin{figure}
\centering
\begin{tabular}{ccc}
\includegraphics[width=0.3\textwidth]{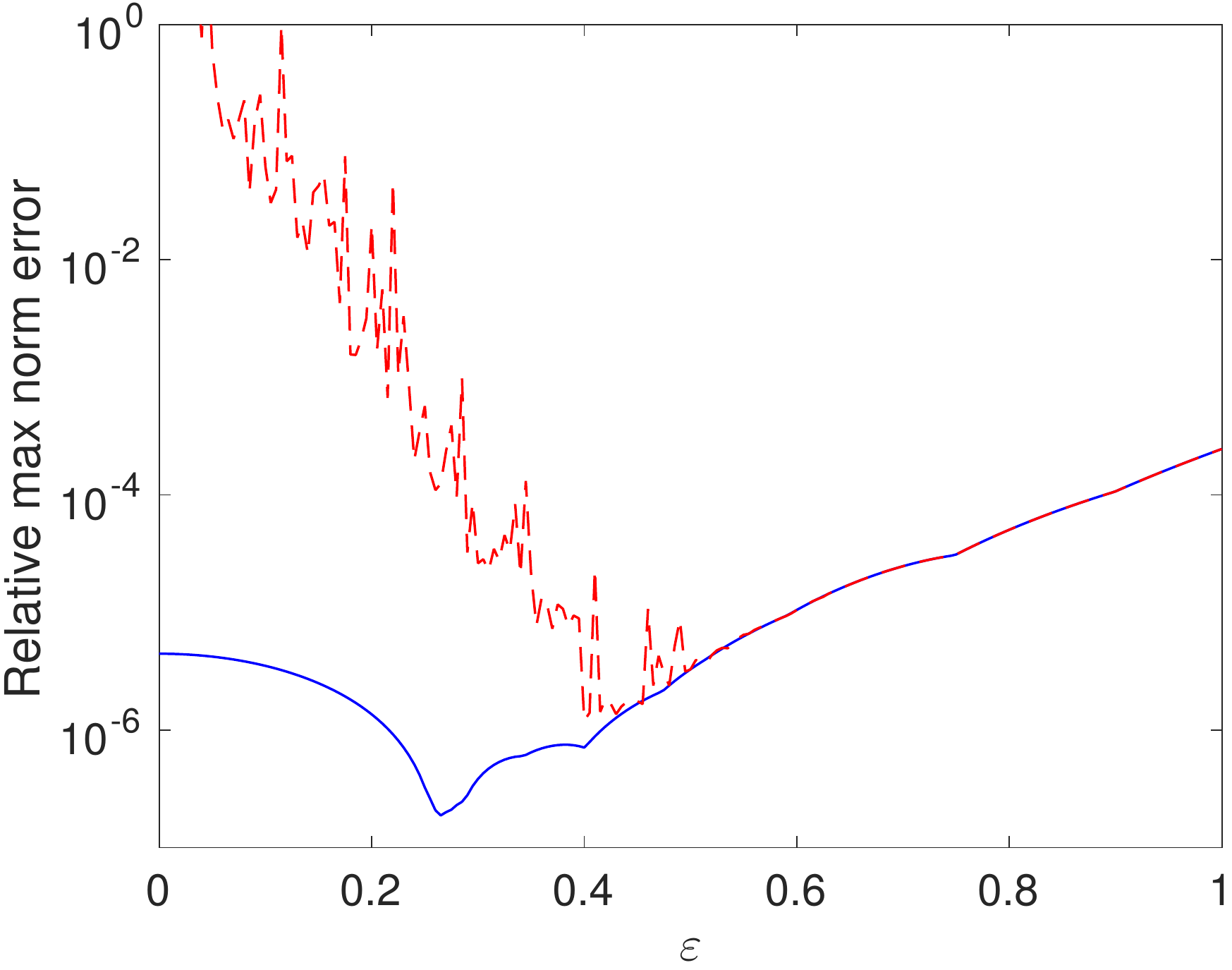} & \includegraphics[width=0.3\textwidth]{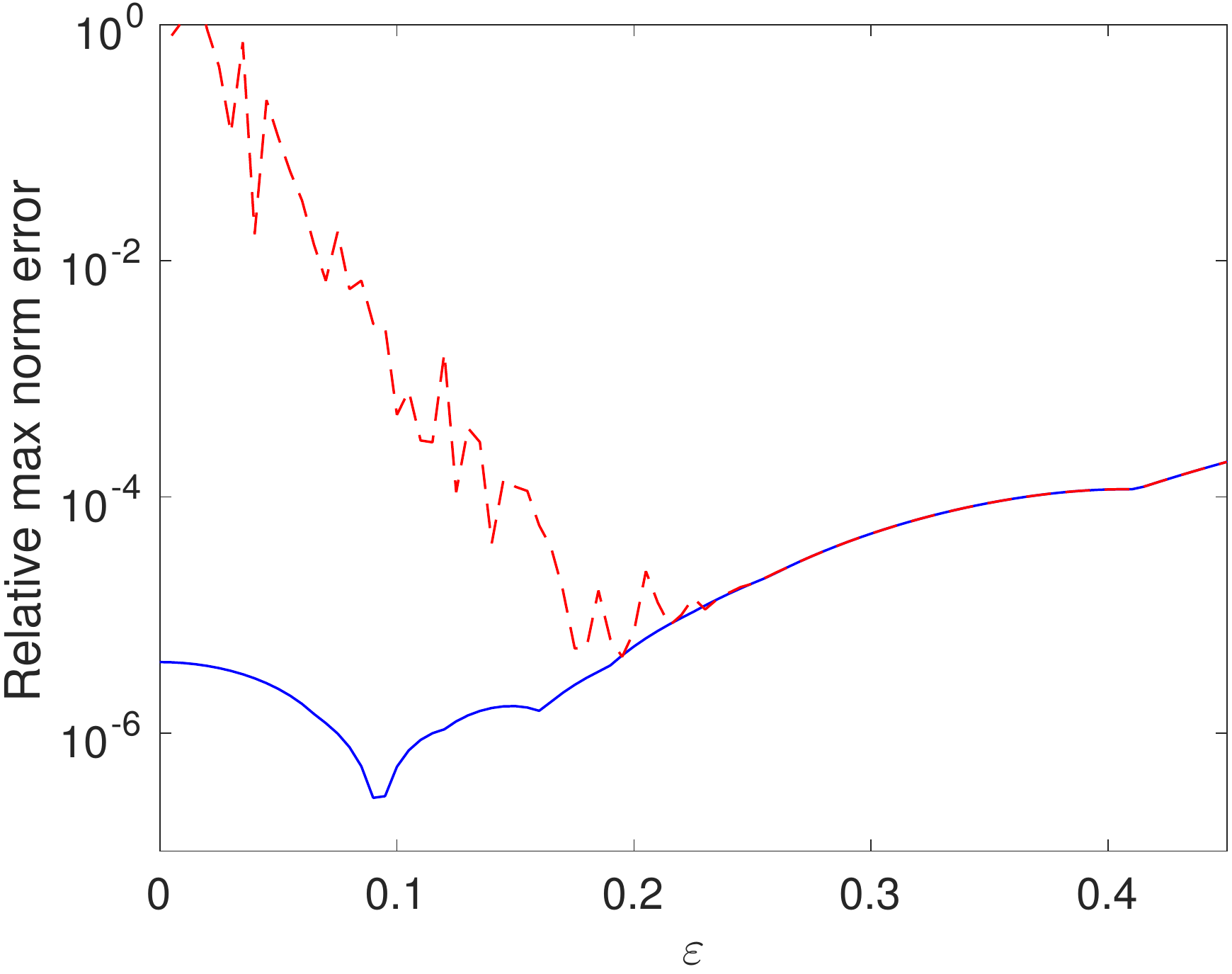} & \includegraphics[width=0.3\textwidth]{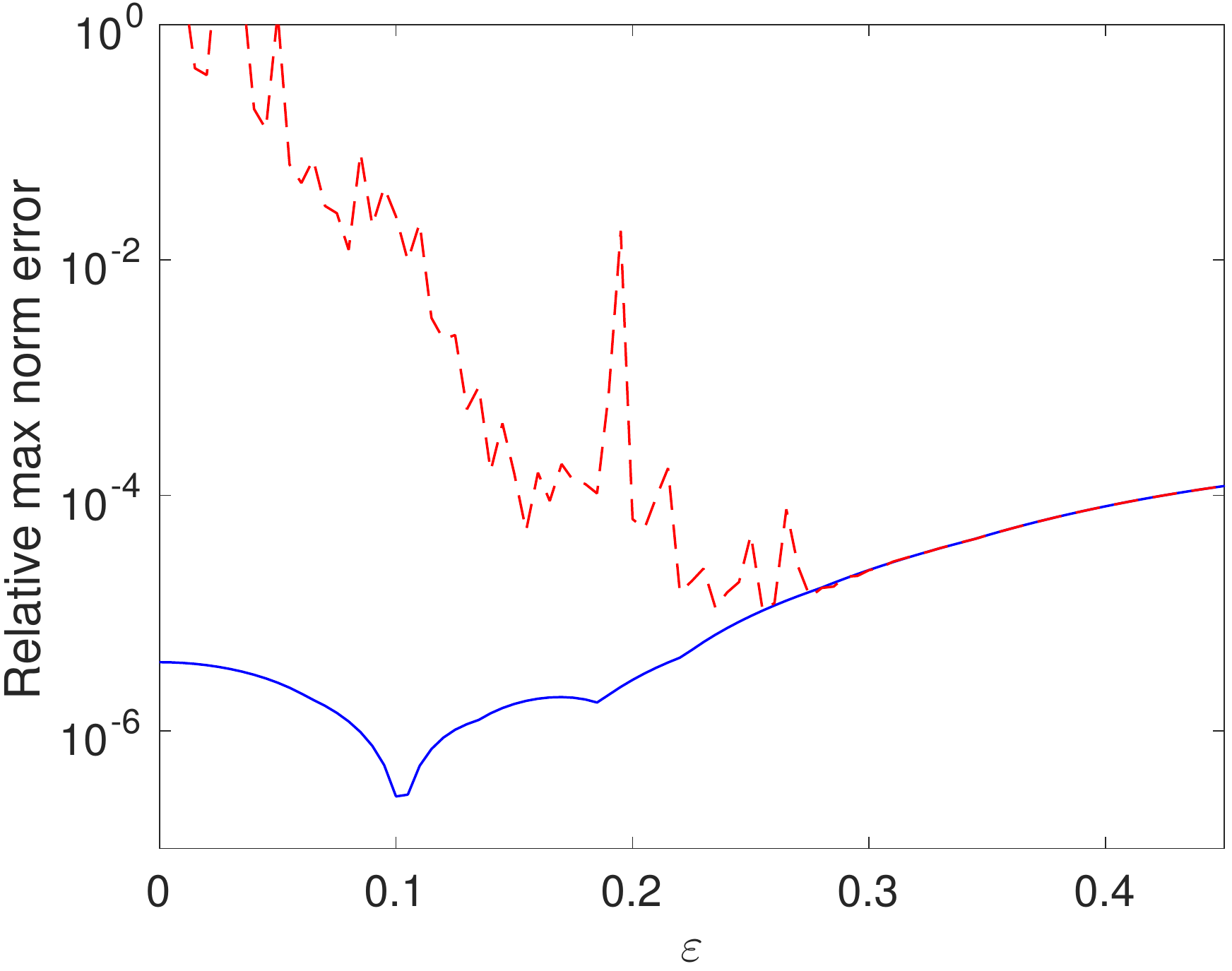} \\
{\small (a) GA} & {\small (b) IQ} & {\small (c) MQ}
\end{tabular}
\caption{Comparison of the relative errors (computed interpolant $-$ target function \eqref{eq:TestFunction}) vs.\ $\ep$ using RBF-Direct (dashed lines) and RBF-RA (solid lines).  Note the different scale on the horizontal axis in (a). \label{fig:err_function}}
\end{figure}
Unless otherwise noted, all numerical results in this subsection are for the $N=62$ nodes over the unit disk shown in Figure \ref{fig:nodes} and for the target function
\begin{align}
g(\vx) = g(x,y) = (1-(x^2 + y^2))\left[ \sin\left(\frac{\pi}{2}(y-0.07)\right) -
                  \frac12\cos\left(\frac{\pi}{2}(x+0.1)\right)\right].
\label{eq:TestFunction}
\end{align}
We use $M=41$ evaluation points $X=\{\vx_j\}_{j=1}^N$ scattered over the unit disk in the RBF-RA algorithm and use these points to measure errors in the interpolant at various values of $\ep$.  

\subsubsection{Errors between interpolant and target function}
In the first numerical experiment, we compare the relative error in the RBF interpolant of the target function \eqref{eq:TestFunction} using RBF-Direct and RBF-RA over a range of different $\ep$.  Specifically, we compute the relative error as 
\begin{align}
\max_{1\leq j \leq M} \left| s(\vx_j,\ep)-g(\vx_j)\right|/\max_{1\leq j \leq M}\left|g(\vx_j)\right|,
\label{eq:err_sg}
\end{align} for $s$ computed with the two different techniques.  The results are shown in Figure \ref{fig:err_function}(a)-(c), for the GA, IQ, and MQ radial kernels, respectively.  We see that RBF-Direct works well for all these kernels until ill-conditioning of the interpolation matrices sets in, while RBF-RA allows the interpolants to be computed in a stable manner right down to $\ep=0$.  The observed pattern of a dip in the RBF-RA error plots is common and is explained in~\cite{FZ07,LarFor05}; it is a feature of the RBF interpolant and does not have to do with ill-conditioning of the RBF-RA method.

\subsubsection{Accuracy vs. $\ep$ and $K$}
In all the remaining subsections of \ref{sec:interp_results} except the last one, we focus on the difference between the computed interpolant and the exact interpolant.  Specifically, we compute
\begin{align}
\max_{1\leq j \leq M} \left| s(\vx_j,\ep)-s_{\text{exact}}(\vx_j,\ep)\right|/\max_{1\leq j \leq M}\left|s_{\text{exact}}(\vx_j,\ep)\right|,
\label{eq:err_f}
\end{align}
where $s$ is the interpolant obtained from either RBF-RA, RBF-CP, or RBF-Direct, and $s_{\text{exact}}$ is the `exact' interpolant, computed using multiprecision arithmetic with 200 digits.  Additionally, for brevity, we limit the presented results to the GA kernel as similar results were observed for other radial kernels.  

\begin{figure}[h]
\centering
\begin{tabular}{cc}
\includegraphics[width=0.45\textwidth]{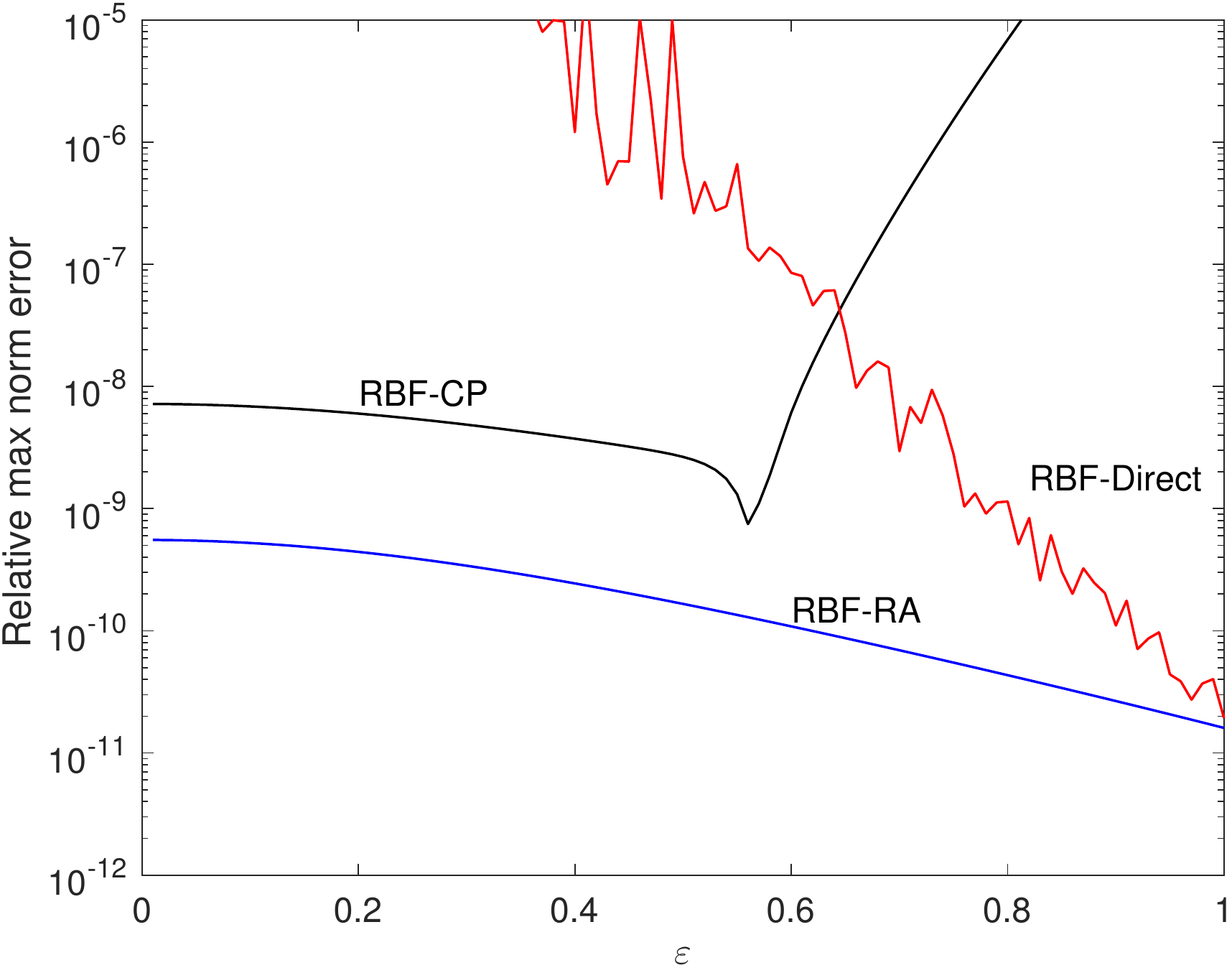} & \includegraphics[width=0.45\textwidth]{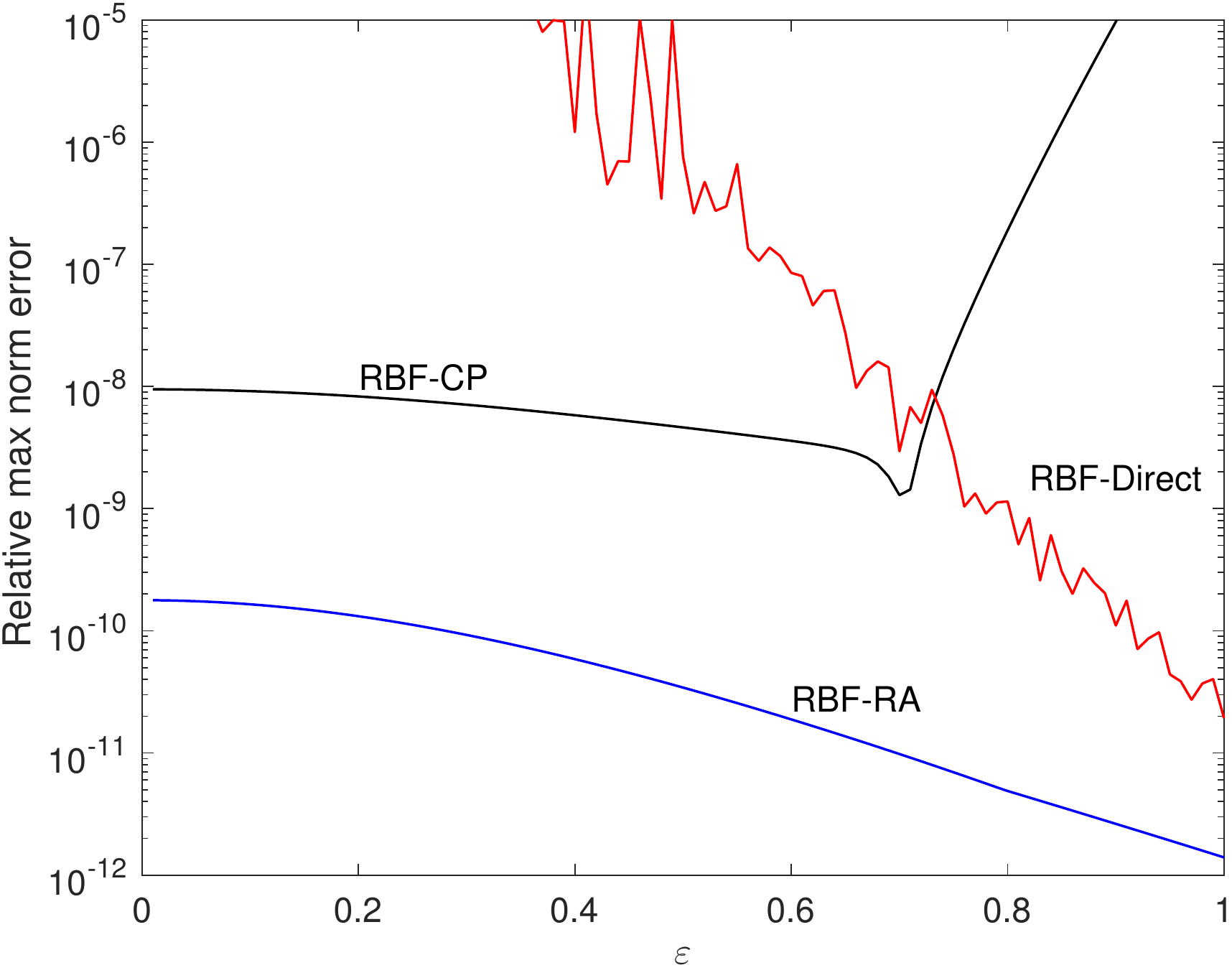} \\
(a) $K/2=22$  & (b) $K/2=32$  \\
\multicolumn{2}{c}{\includegraphics[width=0.45\textwidth]{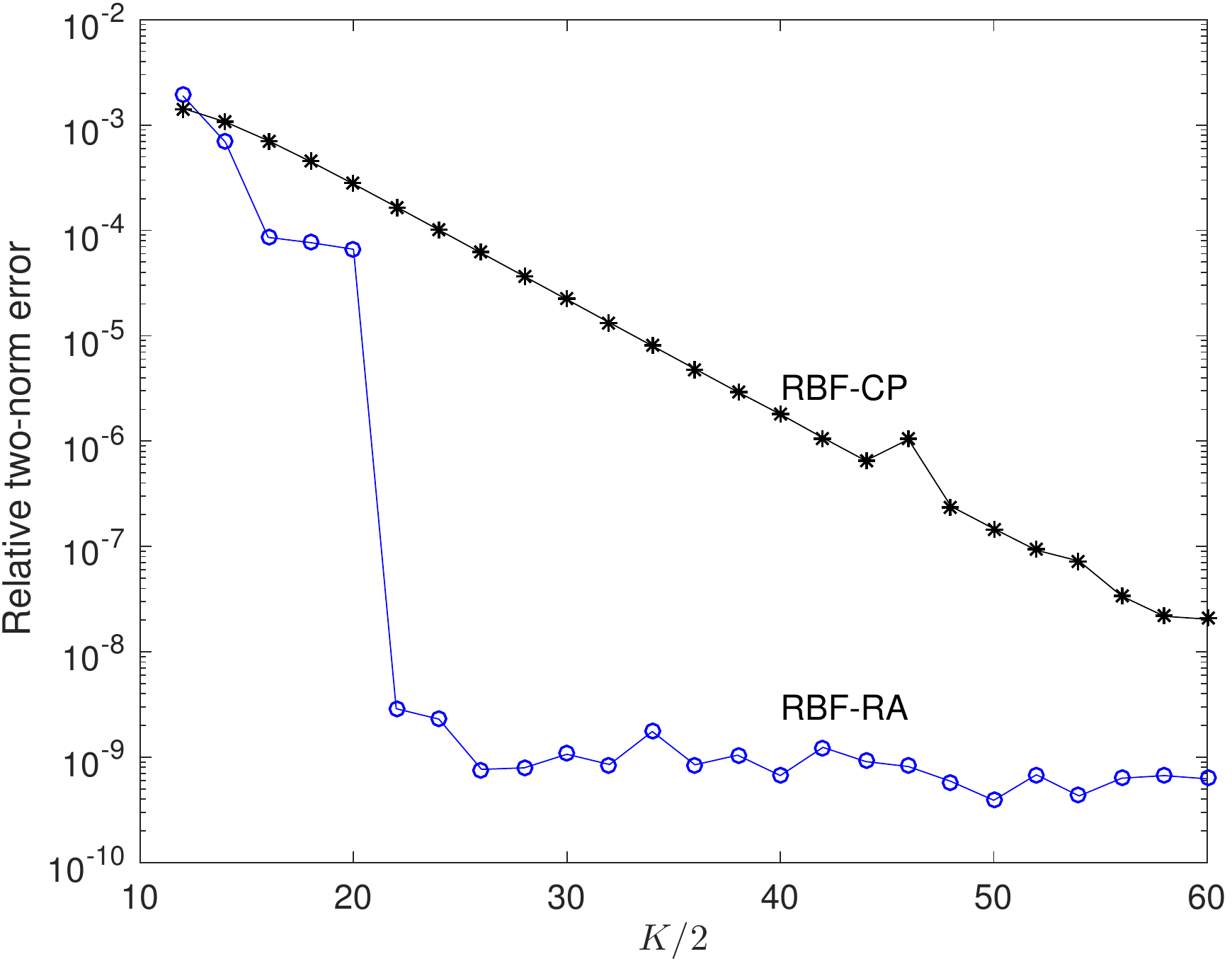}} \\
\multicolumn{2}{c}{(c) Errors over $0 < \ep \leq 0.9$ for various $K$}
\end{tabular}
\caption{Comparison of the errors (computed interpolant $-$ exact interpolant) as the number of evaluations points $K$ on the contour varies.  Figures (a) and (b) show the relative max norm error for $\ep\in[0,1]$ for $K/2=22$ and $K/2=32$ (we use $K/2$ since this is this is the actual number of evaluations of the interpolant that are required).  Figure (c) shows the relative two norm of the error taken only over $0\leq\ep\leq 0.9$ as a function of $K/2$. \label{fig:err_interp}}
\end{figure}

In the first round of experiments, we compare the accuracy of the RBF-RA and RBF-CP algorithms as a function of $\ep$ and $K$ (the number of $\ep$ points on the contour used in both algorithms).  We present the results in terms of $K/2$ since this is the actual total number of evaluations that have to be made of the interpolants (see comment towards the end of Section \ref{sec:vvra}), which is the dominating cost of the algorithm.  Figure \ref{fig:err_interp} (a) shows the relative errors in the interpolants computed using the RBF-RA and RBF-CP methods for $0<\ep\leq 1$ and for $K/2=22$ and $K/2=32$, respectively.  It is immediately obvious from these results that the RBF-RA algorithm gives more accurate results over the entire range of $\ep$ and that this is more prominent as $\ep$ grows.  Additionally, increasing $K$ has only a minor effect on the accuracy of the RBF-RA method, but has a more pronounced effect on the RBF-CP method, with an increased range of $\ep$ for which it is accurate.  We explore the connection between accuracy and $K$ further by plotting in Figure \ref{fig:err_interp} (c) the relative errors in the interpolant for many different $K$.  Here the errors are computed for each $K$, by first computing the max-norm
error \eqref{eq:err_f} and then computing the two-norm of these errors over $0 < \ep \leq 0.9$. We see that the RBF-RA method converges rapidly with $K$, while the RBF-CP method converges much slower (but still at a geometric rate).  The reason the error stops decreasing in the RBF-RA method around $10^{-9}$ is that this is about the relative accuracy that can be achieved using RBF-Direct to compute the interpolants on the contour in the complex $\ep$-plane.

\subsubsection{RBF-RA vs. Multiprecision RBF-Direct}
\begin{figure}
\centering
\begin{minipage}{0.48\textwidth}
\centering
\includegraphics[width=0.99\textwidth]{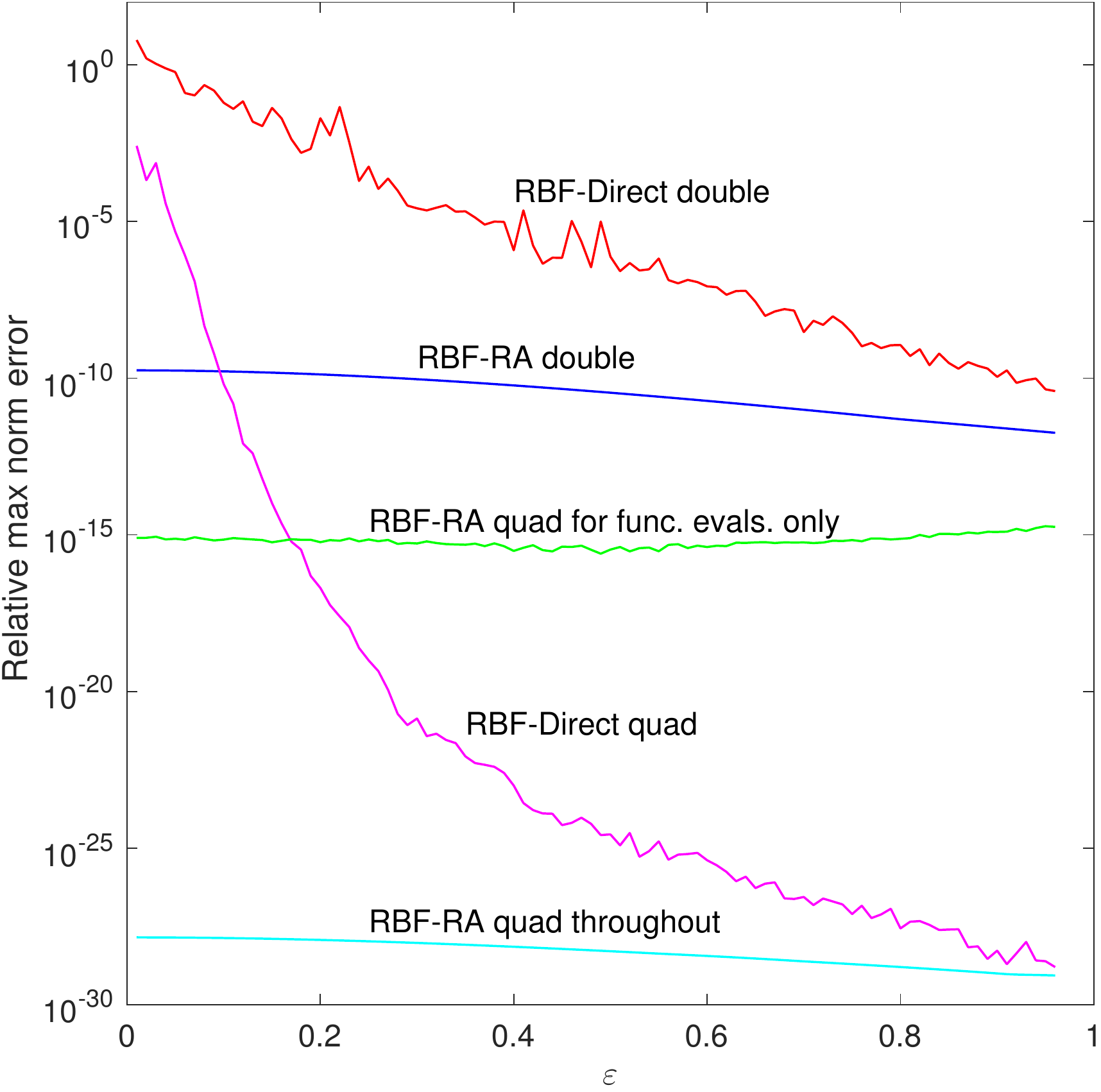}\vspace{1pt}
\caption*{(a)}
\end{minipage}
\begin{minipage}{0.47\textwidth}
\centering
\includegraphics[width=0.99\textwidth]{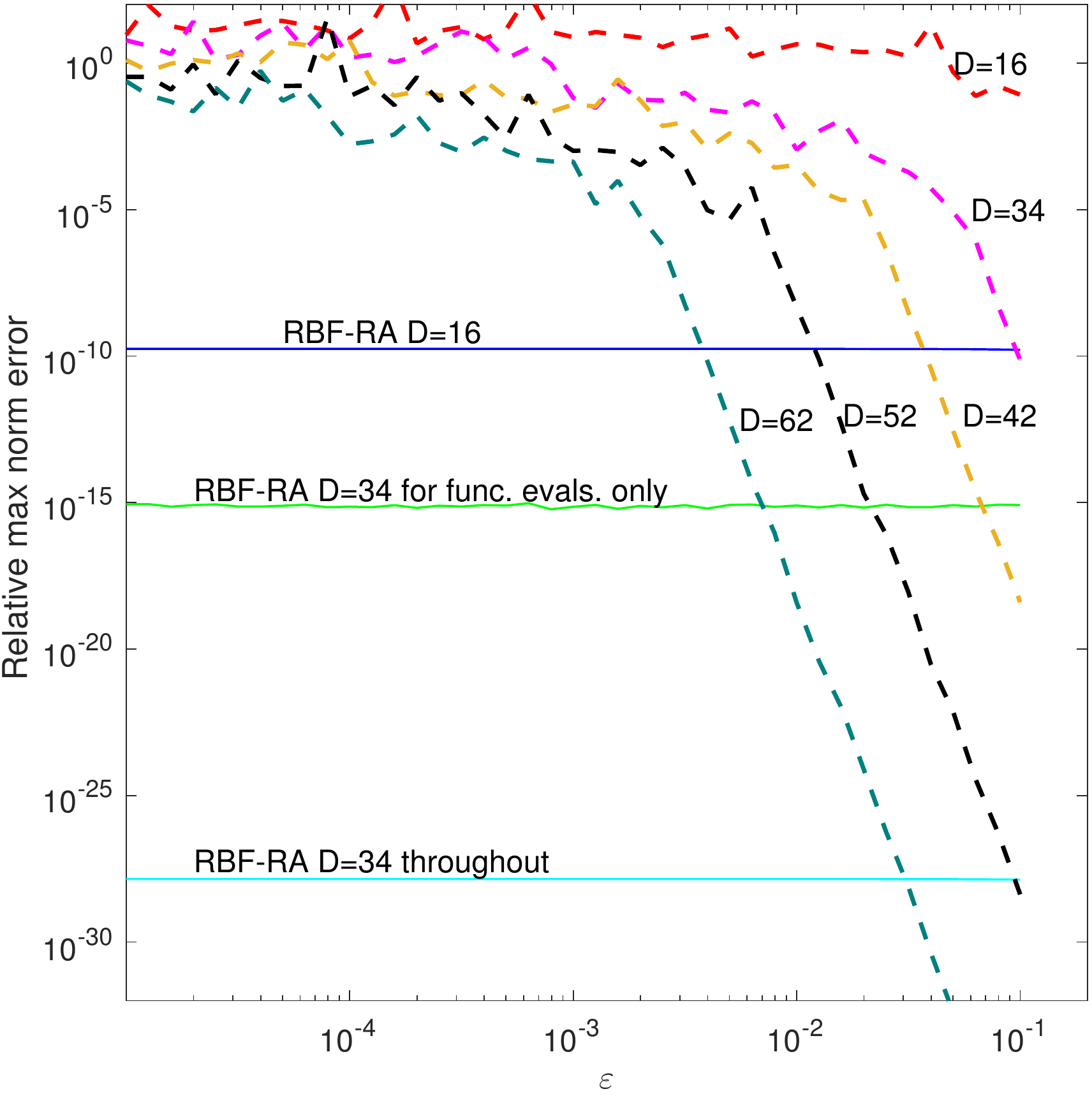}
\caption*{(b)}
\end{minipage}
\vspace{-10pt}
\caption{(a) Comparison of the errors (computed interpolant $-$ exact interpolant) using double precision and quad precision arithmetic.  The top RBF-RA quad curve corresponds to using quad precision only when evaluating the interpolant, whereas the bottom quad curve corresponds to using quad precision for all the computations in the RBF-RA method.  (b) Similar to (a) but for multiple precision arithmetic using D digits (here D=16 is double precision and D=34 is quad precision).  Here RBF-Direct results are in the dashed lines; note also the logarithmic scale in $\ep$. \label{fig:err_mp}}
\end{figure}

An all too common way to deal with the ill-conditioning associated with flat RBFs is to use multiprecision floating point arithmetic with RBF-Direct.  While this does allow one to consider a larger range of $\ep$ in applications, it comes at a higher computational cost as the computations have to be done using special software packages instead of directly in hardware.  Typically, quad precision (approximately 34 decimal digits) is used since floating point operations on these numbers can combine two double precision numbers, which can improve the cost.  Figure \ref{fig:err_mp} (a) compares the errors in the interpolants computed using RBF-RA and RBF-Direct with both double and quad precision arithmetic, with the latter being computed using the Advanpix multiprecision \matlab toolbox~\cite{Advanpix}.  We see from the figure that the range of $\ep$ that can be considered with RBF-Direct and quad precision increases, but that the method is still not able to consider the full range.  There is nothing to prevent quad precision from also being used with RBF-RA and in Figure \ref{fig:err_mp} (a) we have included two results with quad precision.  The first uses quad precision only to evaluate the interpolant (which is the step that has the potential to be ill-conditioned), with all subsequent computations of the RBF-RA method done in double precision.  The second uses quad precision throughout the entire RBF-RA algorithm.  We see from the results that the errors in the interpolant for both cases are now much lower than the double precision case.  In Figure \ref{fig:err_mp} (b) we further explore the use of multiprecision arithmetic with the RBF-Direct algorithm (again using Advanpix) by looking at the error in the interpolant as the number of digits used is increased.  Now, we  consider $10^{-5} < \ep \leq 10^{-1}$ to clearly illustrate that the ill-conditioning of RBF-Direct cannot be completely overcome with this technique, but that it can be completely overcome with RBF-RA.

\subsubsection{Contours passing close to poles}
\begin{figure}
\centering
\begin{tabular}{cc}
\includegraphics[width=0.42\textwidth]{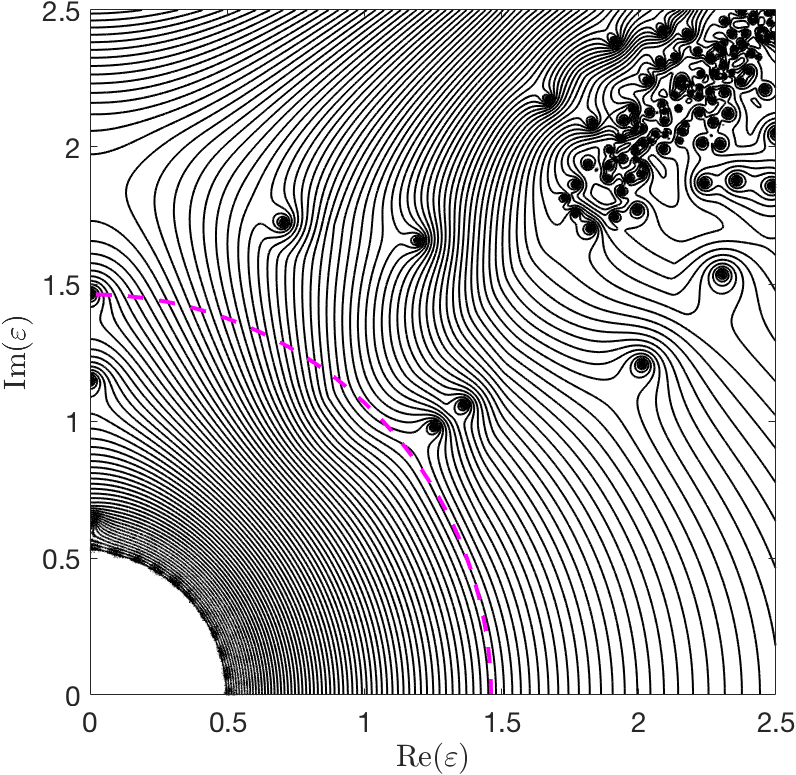} & \includegraphics[width=0.52\textwidth]{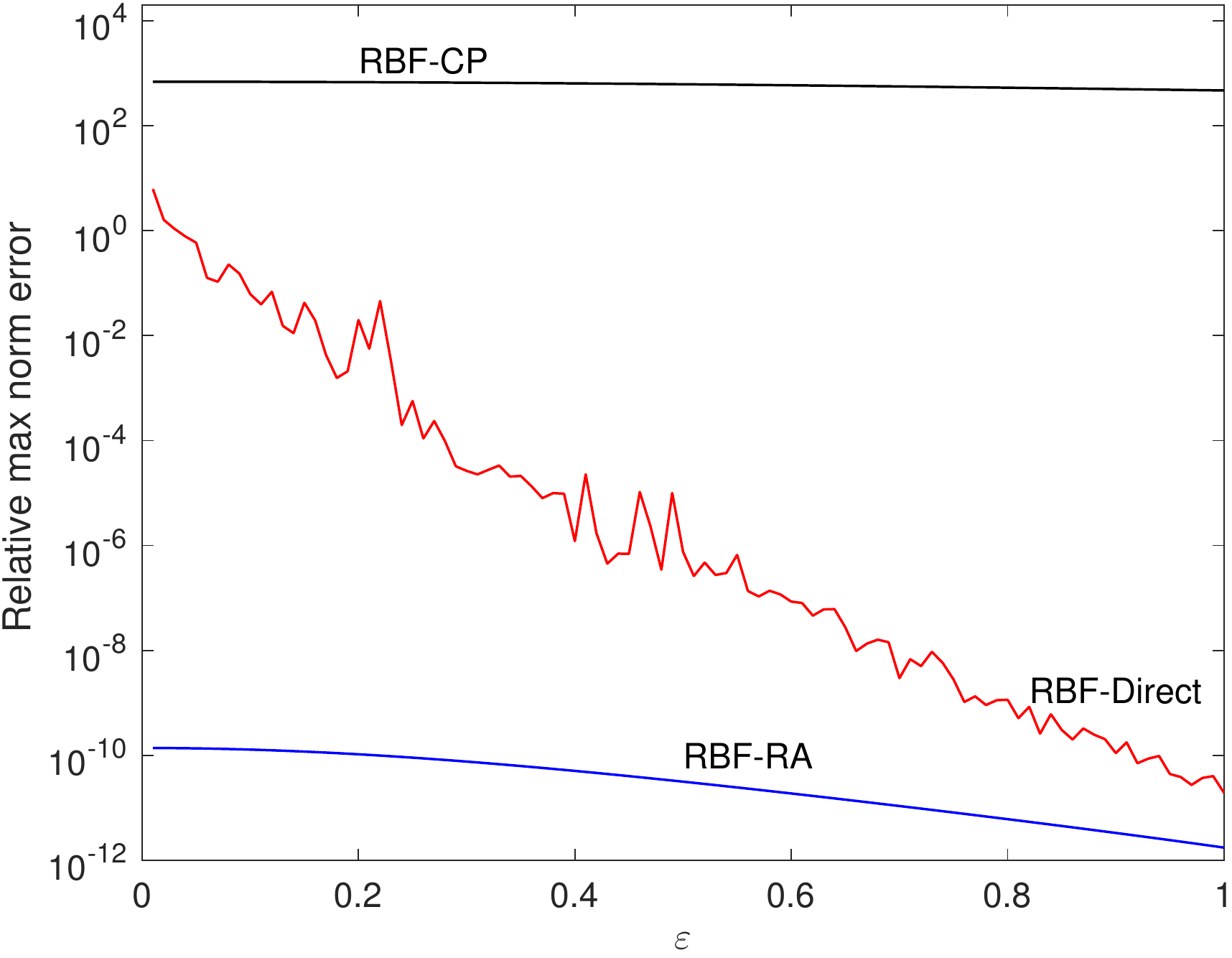} \\
(a)  & (b) 
\end{tabular}
\caption{(a) Contour plot of $\log_{10}(\cond(A(\ep)))$ (similar to Figure \ref{fig:cond_surface} (b)) with an evaluation contour (dashed line) running very close to a pole on the imaginary axis.  (b) Comparison of the errors (computed interpolant $-$ exact interpolant)  associated with the RBF-RA and RBF-CP algorithms using the evaluation contour in part (a). Here we used $K/2=32$ evaluation points on the contour; doubling this number does not appear to improve the RBF-CP results at all. \label{fig:err_pole}}
\end{figure}

To demonstrate the robustness of the new RBF-RA algorithm we consider a case where the evaluation contour runs very close to a pole in the complex $\ep$-plane.  For the $N=62$ node example we are considering in these numerical experiments, there is a pole at $\ep\approx 1.4617904771448i$.  We choose for both the RBF-RA and RBF-CP algorithms a circular evaluation contour centered at the origin of radius $\ep = 1.4618$. This contour is superimposed on a contour plot of the condition number in Figure \ref{fig:err_pole} (a) (see the dashed curve).  Figure \ref{fig:err_pole} (b) shows the max-norm errors (again together with RBF-Direct for comparison) in the interpolants \eqref{eq:err_f} computed with both algorithms using this contour.  We see that the RBF-CP algorithm gives entirely useless results for this contour, while the RBF-RA algorithm performs as good (if not slightly better) than the case where the contour does not run close to any poles.

\subsubsection{Vector valued rational approximation vs.\ rational interpolation}
\begin{figure}
\centering
\includegraphics[width=0.45\textwidth]{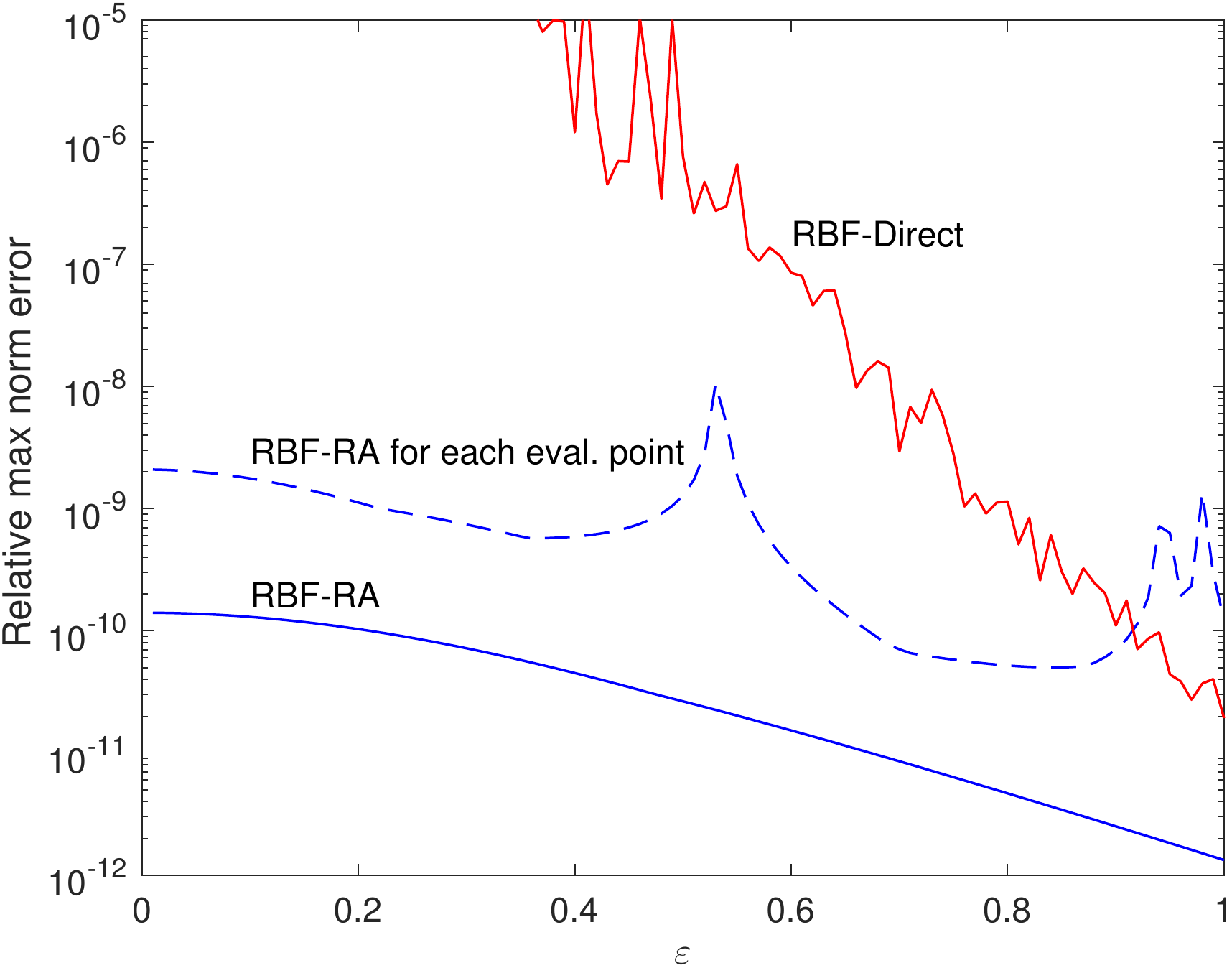}
\caption{Comparison of the errors (computed interpolant $-$ exact interpolant) when using RBF-RA with all the evaluation points (solid line marked RBF-RA), when applying RBF-RA separately to each evaluation point (dashed line), and when using RBF-Direct.  In the second case, a rational approximant is computed separately for each entry of the vector-valued function.  In both RBF-RA  cases, we set $K/2 = 40$. \label{fig:err_ra_ri}}
\end{figure}

A seemingly simpler approach to obtain vector valued rational approximations to the interpolant $\vs(\ep)$ in \eqref{eq:rbf_mero} is to use rational interpolation (or approximation) of each of the $M$ entries of $\vs(\ep)$ \textit{separately}~\cite{GonnetPachonTrefethen}, instead of the RBF-RA procedure that couples the entries together.  However, we have found that this does not produce as accurate results and can lead to issues with the approximants.  Figure \ref{fig:err_ra_ri} illustrates this by comparing the errors that result from approximating $\vs$ using the standard RBF-RA procedure to that of using RBF-RA separately for each of the $M$ entries of $\vs(\ep)$ (which amounts to computing a rational interpolant of each entry).  We first see from this figure that the RBF-RA procedure is at least an order of magnitude more accurate than the rational interpolation approach.  Second, we see a few values of $\ep$ where the error spikes in the rational interpolation method.  These spikes correspond to spurious poles, or ``Froissart doublets''~\cite{Froissart}, appearing near or on the real $\ep$-axis.  Froissart doublets are especially common in rational approximation methods where the input contains noise~\cite{Beckermann201591}, which occurs in our application at roughly the unit roundoff of the machine times the condition number of the RBF interpolation matrix on the evaluation contour of the RBF-RA algorithm.  The least squares nature of determining a common denominator in the RBF-RA method appears to significantly reduce the presence of these spurious poles.  In fact, we have not observed the presence of Froissart doublets in any of our experiments with RBF-RA. 

\subsubsection{Timing example for a 3D problem}
\begin{figure}
\centering
\begin{minipage}{0.48\textwidth}
\centering
\includegraphics[width=0.99\textwidth]{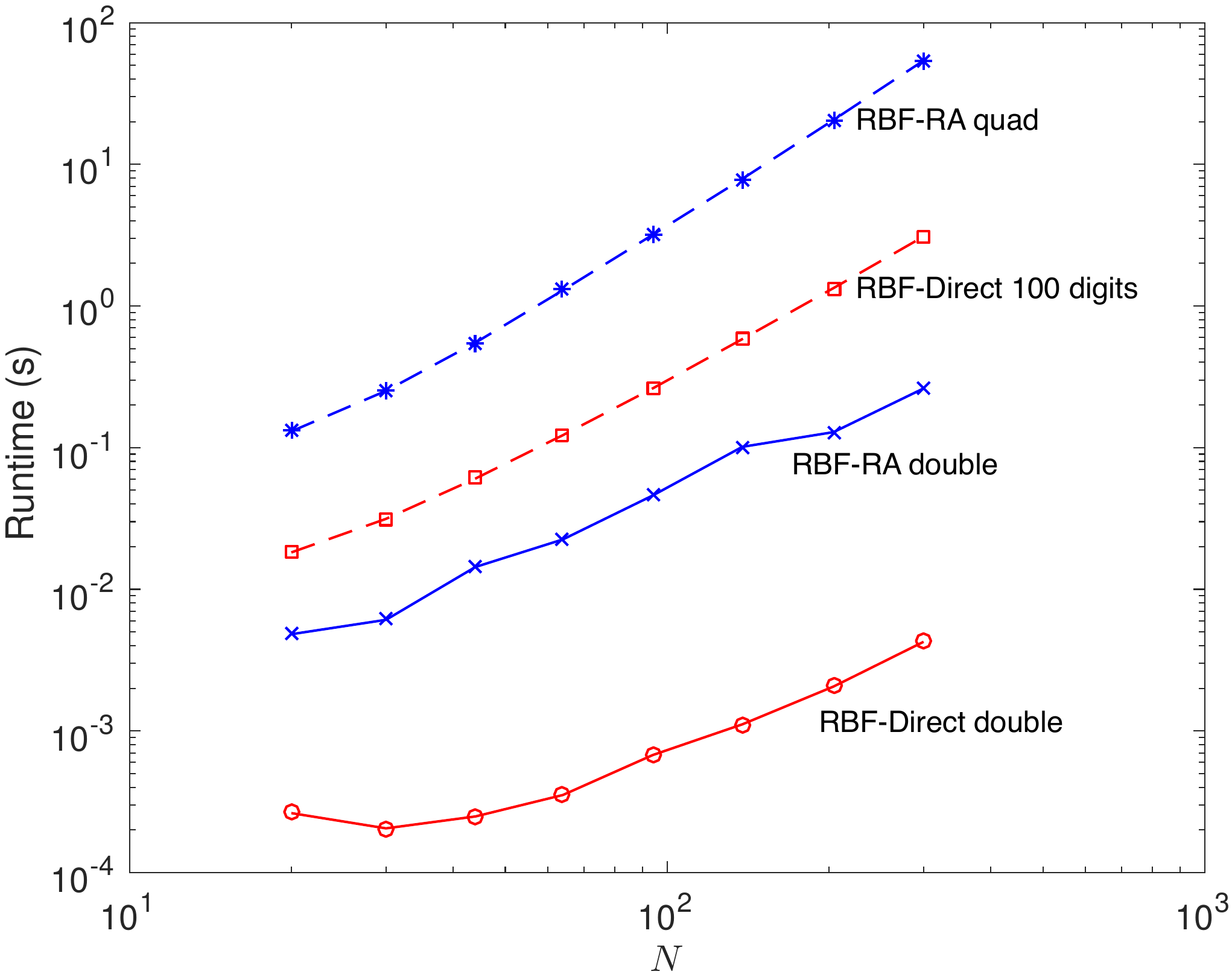}\vspace{1pt}
\caption*{(a) One value of $\ep$.}
\end{minipage}
\begin{minipage}{0.481\textwidth}
\centering
\includegraphics[width=0.995\textwidth]{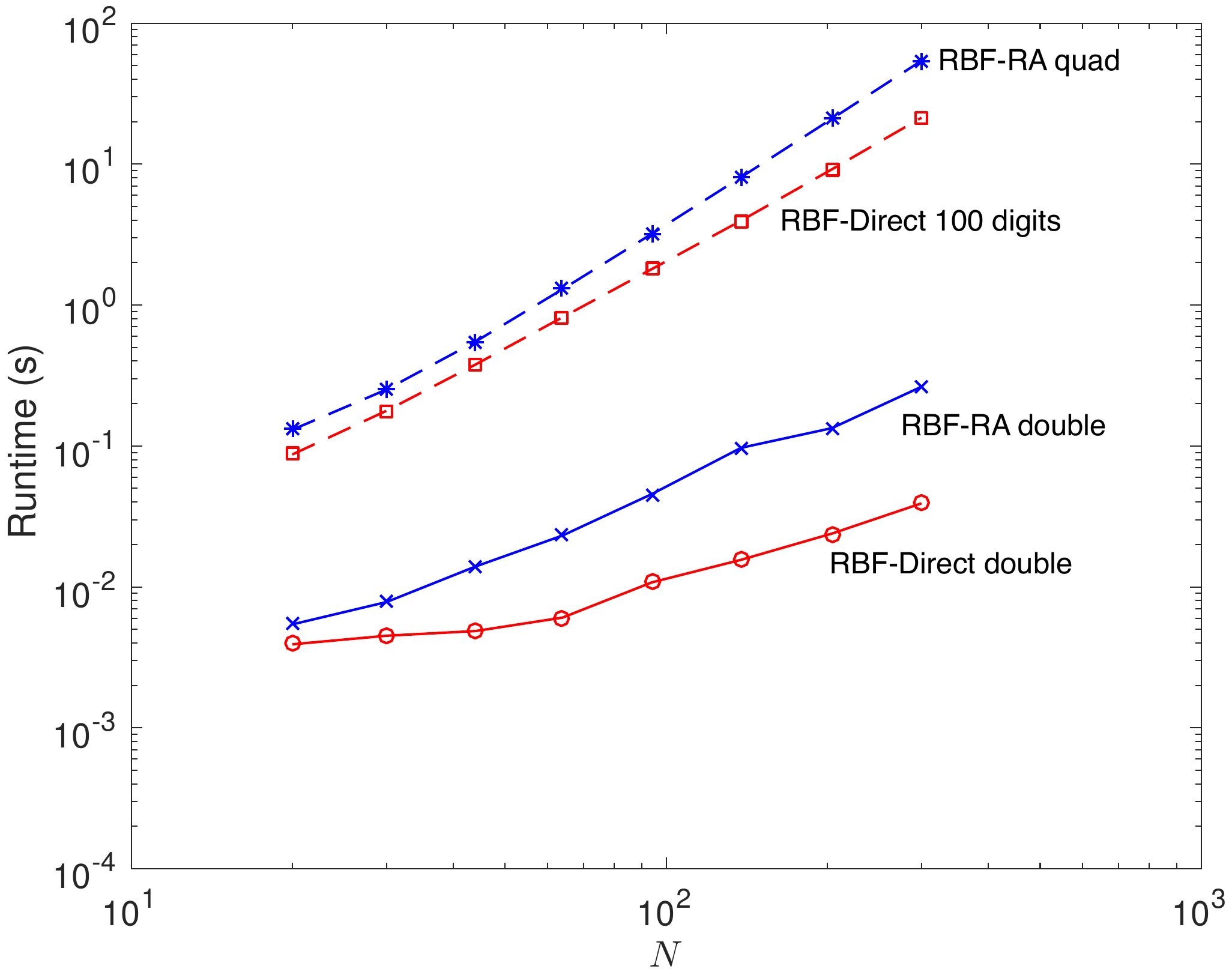}
\caption*{(b) Ten values of $\ep$. }
\end{minipage}
\vspace{-10pt}
\caption{Measured wall clock time (in seconds) as a function of the number of 3-D Halton nodes $N$ for the RBF-RA (with $K/2=32$) and RBF-Direct methods.  (a) A comparison of the methods for computing the interpolant at $M=N$ evaluation points and for one value of $\ep$ (set to $\ep=10^{-2}$).  (b) Same as (a), but for computing the interpolant at ten values of $\ep$ (set to $\ep_j = 10^{-j}$, $j=0,\ldots,9$).  The dashed lines were computed using multiprecision arithmetic, with $D=100$ digits for RBF-Direct and $D=34$ (quad precision) for RBF-RA.  All timings were done using \matlab 2016a on a 2014 MacBook Pro with 16 GB of RAM and an 3GHz Intel Core i7 processor without explicit parallelization.  The multiprecision computations were done using the Advanpix multiprecision \matlab toolbox~\cite{Advanpix}.\label{fig:timings}}
\end{figure}

Finally, we compare the computational cost of the RBF-RA method to that of the RBF-Direct method. For this comparison, we use Halton node sets over the unit cube in $\mathbb{R}^3$ of increasing size $N$ (generated with the \matlab function \texttt{haltonset}).  For the evaluation points, we also use Halton node sets with the same cardinality as the nodes (i.e.\ $M=N$),  but shift them so they don't agree with the nodes; having $M=N$ is a common scenario in generating RBF-FD/HFD formulas.  The target function $g$ is not important in this experiment, as we are only concerned with timings.  Figure \ref{fig:timings} (a) shows the measured wall clock times of two methods for evaluating the interpolant at $\ep=10^{-2}$.  Included in this figure are the wall clock times also for computing the interpolants using multiprecision arithmetic, with $D=100$ digits for RBF-Direct and with $D=34$ (quad precision) for RBF-RA.  For the $N>100$ and $\ep=0.01$, it is necessary to switch to multiprecision arithmetic with RBF-Direct to get a meaningful result, whereas RBF-RA in double precision has no issues with ill-conditioning.  While $D=100$ is larger than is necessary for RBF-Direct with $\ep=10^{-2}$, it is reasonable for smaller $\ep$ (see Figure \ref{fig:err_mp} (b)), also the timings do not go down very much by decreasing $D$.  The quad precision results are included for comparison purposes with RBF-Direct in multiprecision mode; they are not necessary to obtain an accurate result for these values of $N$.  For the double precision computations, we see that the cost of RBF-RA is about 100 times that of RBF-Direct.  However, for small $\ep$, the comparison should really be made between RBF-Direct using multiprecision arithmetic, in which case RBF-RA is about an order of magnitude more efficient.  Based on the timing results in~\cite[Section 6.3]{FoLePo13}, the computational cost of the RBF-RA method is about a factor of two or three larger than the RBF-QR method and about a factor of ten larger than the RBF-GA method, which (at about 10 times the cost of RBF-Direct) is the fastest of the stable algorithms (for GA RBFs).  

One added benefit of the RBF-RA method is that evaluating the interpolant at multiple values of $\ep$ (which is required for some shape parameter selection algorithms) comes at virtually no additional computational cost.  The same is not true for the RBF-Direct method.  We demonstrate this feature in Figure \ref{fig:timings} (b), where the wall clock times of the algorithms are displayed for evaluating the interpolants at $\ep_j = 10^{-j}$, $j=0,\ldots,9$.   

The dominant cost of the RBF-RA method comes from evaluating the interpolant on the contour, which can be done in parallel.  The timings presented here have made no explicit use of parallelization for these evaluations, so further improvements to the efficiency the method over RBF-Direct are still possible.

\subsection{RBF-HFD results}
\begin{figure} 
\centering
\begin{minipage}{.45\textwidth}
\begin{overpic}[width=0.9\textwidth]{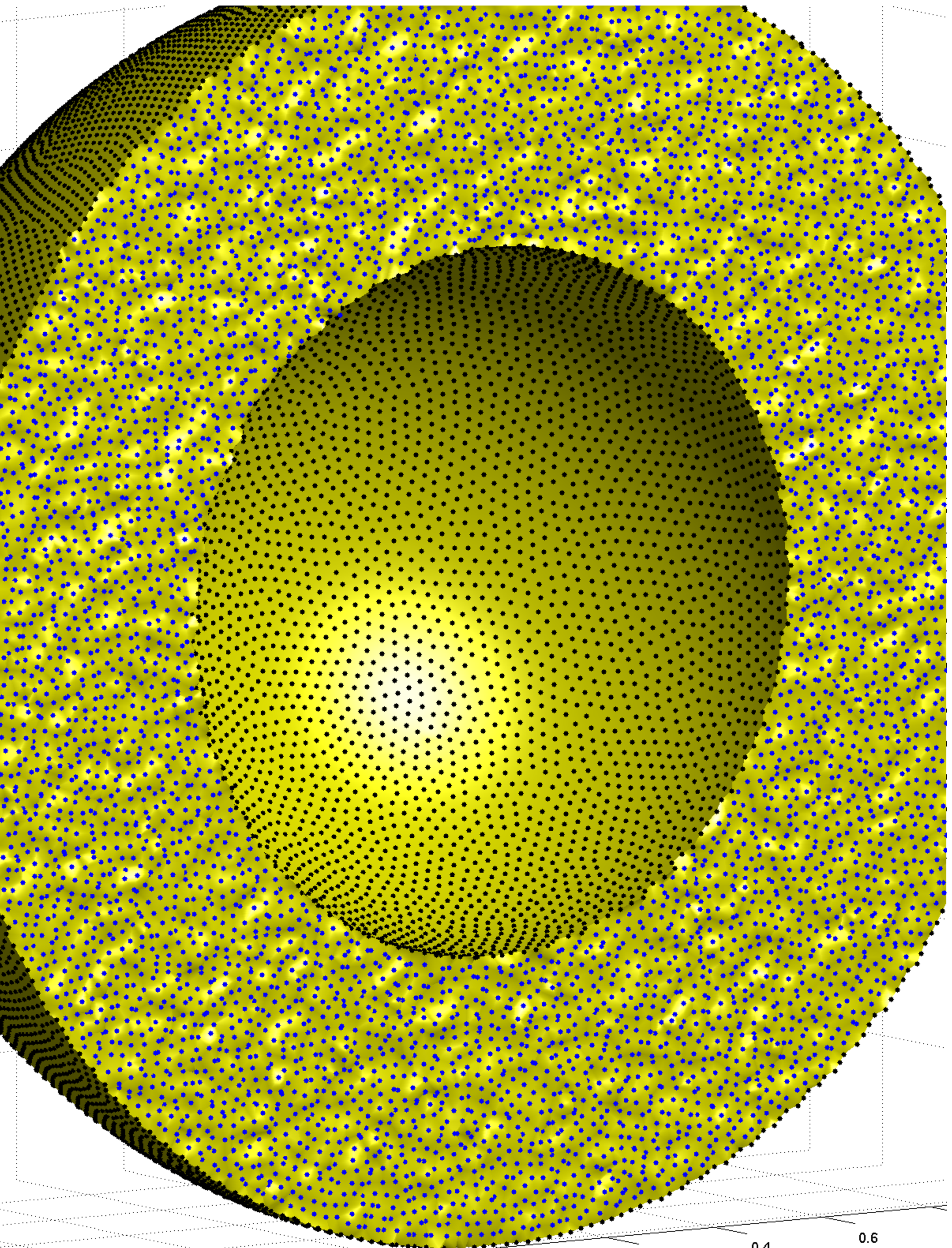}
\end{overpic}
\end{minipage}
\begin{minipage}{0.45\textwidth}
\begin{overpic}[width=0.6\textwidth]{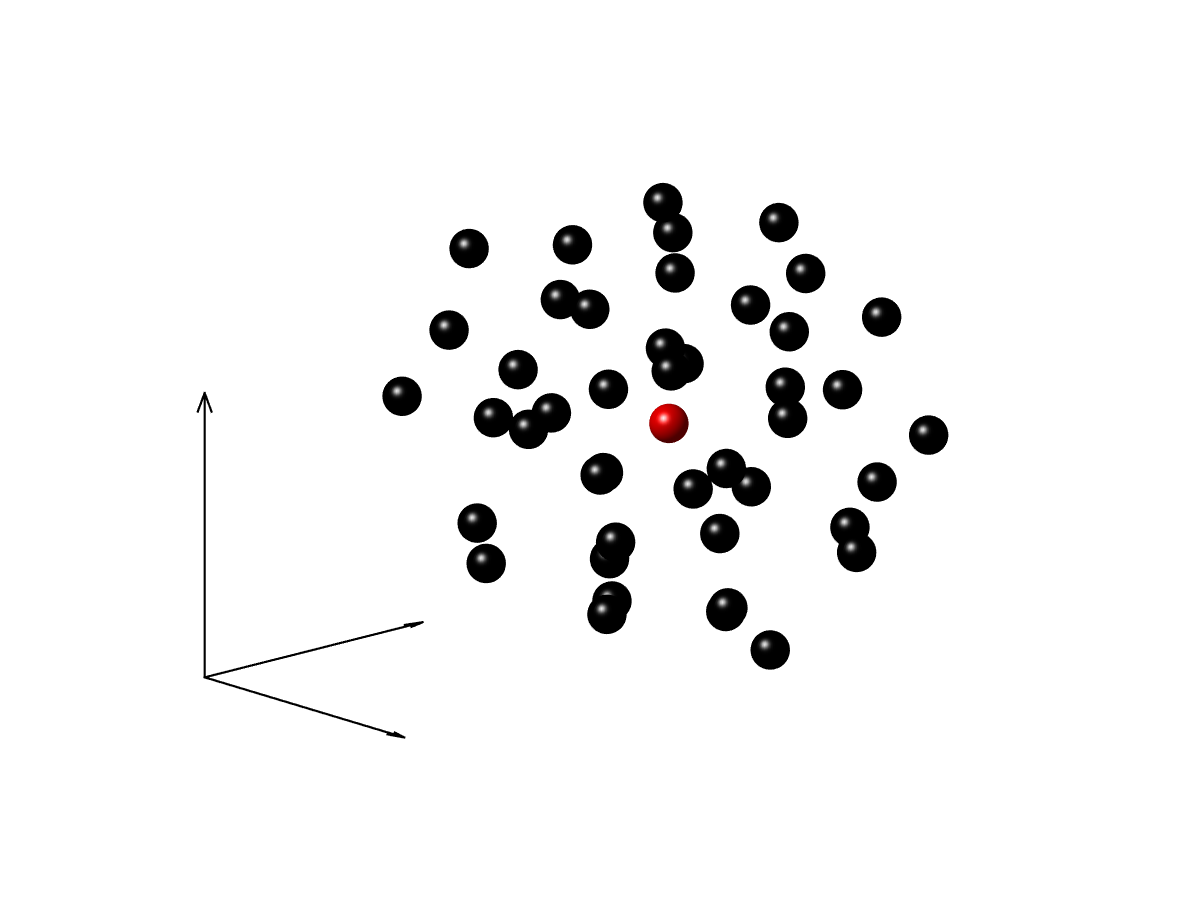}
\put(-35,90){\rotatebox{-32}{\Large$\xrightarrow{\text{\phantom{foooooooooooooooooo}}}$}}
\end{overpic}
\end{minipage}
\caption{Illustration of the global and local node set used in the RBF-HFD method for solving Poisson's equation in a spherical shell with radius $0.55\leq r \leq 1$, which mimics the aspect radius of the Earth's mantle.  The left figure shows the shell split open, with part of the  $5\cdot 10^5$ global node set marked by small solid spheres.  The right shows a small subset of the nodes, the HFD stencil, with the `center' node marked in red.\label{fig:rbf_fd_illustration}}
\end{figure}

In this section we consider the application of the RBF-RA method to computing RBF-HFD weights for the 3-D Laplacian and use these to solve Poisson's equation in a spherical shell.  Specifically, we are interested solving
\begin{align}
\Delta u &= g,\; \text{in}\; \Omega = \left\{(x,y,z)\in\mathbb{R}^d \Bigl | 0.55 \leq \sqrt{x^2 + y^2 + z^2} \leq 1\right\},
\label{eq:Poisson}
\end{align}
subject to Dirichlet boundary conditions on the inner and outer surfaces of the shell.  Here we take the exact solution to be
\begin{align*}
u(\lambda,\theta,r) = \sin\left(\frac{20\pi}{9}\left(r-\frac{11}{20}\right)\right)\left[Y_{6}^{0}(\lambda,\theta) + \frac{14}{11}Y_{6}^{5}(\lambda,\theta)\right],
\end{align*}
where $(\lambda,\theta,r)$ are spherical coordinates and $Y_{\ell}^m$ denote real-valued spherical harmonics of degree $\ell$ and order $m$.  This $u$ is used to generate $g$ in \eqref{eq:Poisson} to set up the problem.  We use $500,000$ global nodes to discretize the shell\footnote{These nodes were generated by Prof.\ Douglas Hardin at Vanderbilt University using a modified version of the method discussed in\cite{Borodachov2014}.}; see the left image of Figure \ref{fig:rbf_fd_illustration} for an illustration of the nodes.  We denote the nodes interior to the shell by $\Xi^{\text{int}} = \{\boldsymbol{\xi}_k\}_{k=1}^P$ and the nodes on the boundary by $\Xi^{\text{bnd}} = \{\boldsymbol{\xi}_k\}_{k=P+1}^{P+Q}$.  For our test problem, $P=453,405$ and $Q=46,595$.

The procedure for generating the RBF-HFD formulas is as follows: For $k=1,\ldots,P$ repeat the following
\begin{enumerate}
\item Select the $N-1$ nearest neighbors of $\boldsymbol{\xi}_k$ from $\Xi^{\text{int}} \cup \Xi^{\text{bnd}}$ (using a KD-tree for efficiency), with $N << P$.  These nodes plus $\boldsymbol{\xi}_k$ form the set $\Xh=\{\vxc_i\}_{i=1}^N$ in \eqref{eq:rbf_hfd}, with the convention that $\vxc_1=\boldsymbol{\xi}_k$.
\item From the set $\Xh$, select the $L < N$ nearest-neighbors to $\vxc_1$.  These nodes form the set $\Yh=\{\vyc_j\}_{j=1}^L$ in \eqref{eq:rbf_hfd}.
\item Use $\Xh$, $\Yh$, and  $\calL = \Delta$ in \eqref{eq:rbf_hfd_mero}, then apply the vector-valued rational approximation algorithm to compute the weights for given values of $\ep$.  Note that standard Cartesian coordinates can be used in applying $\calL$ in \eqref{eq:rbf_hfd_mero}.
\end{enumerate}
In the numerical experiments, we set $N=45$ and $L=20$ and use the IQ radial kernel.  We also set $K/2=64$ and choose $n=K/4$.  

\begin{figure}
\centering
\begin{tabular}{cc}
\includegraphics[width=0.47\textwidth]{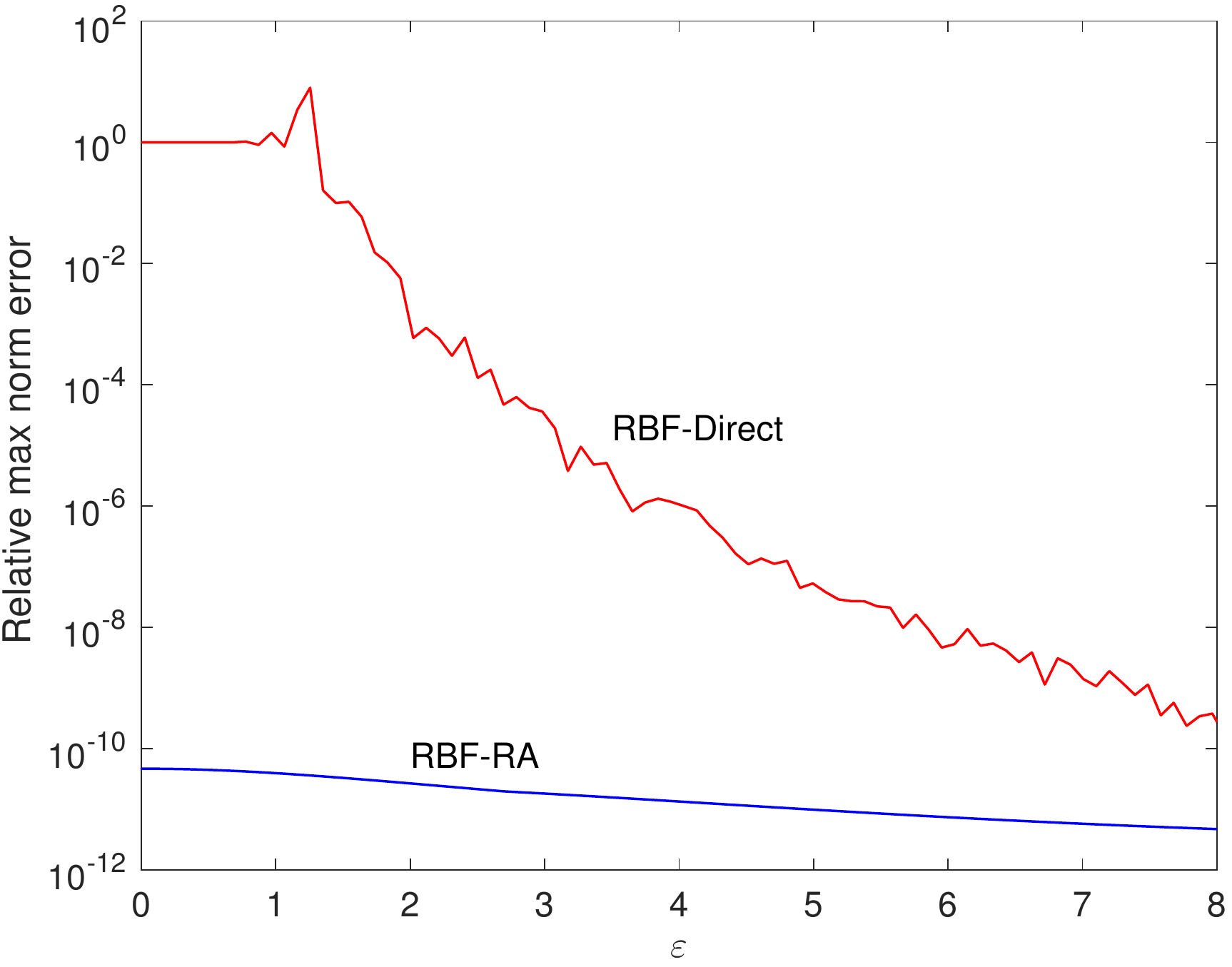} & \includegraphics[width=0.46\textwidth]{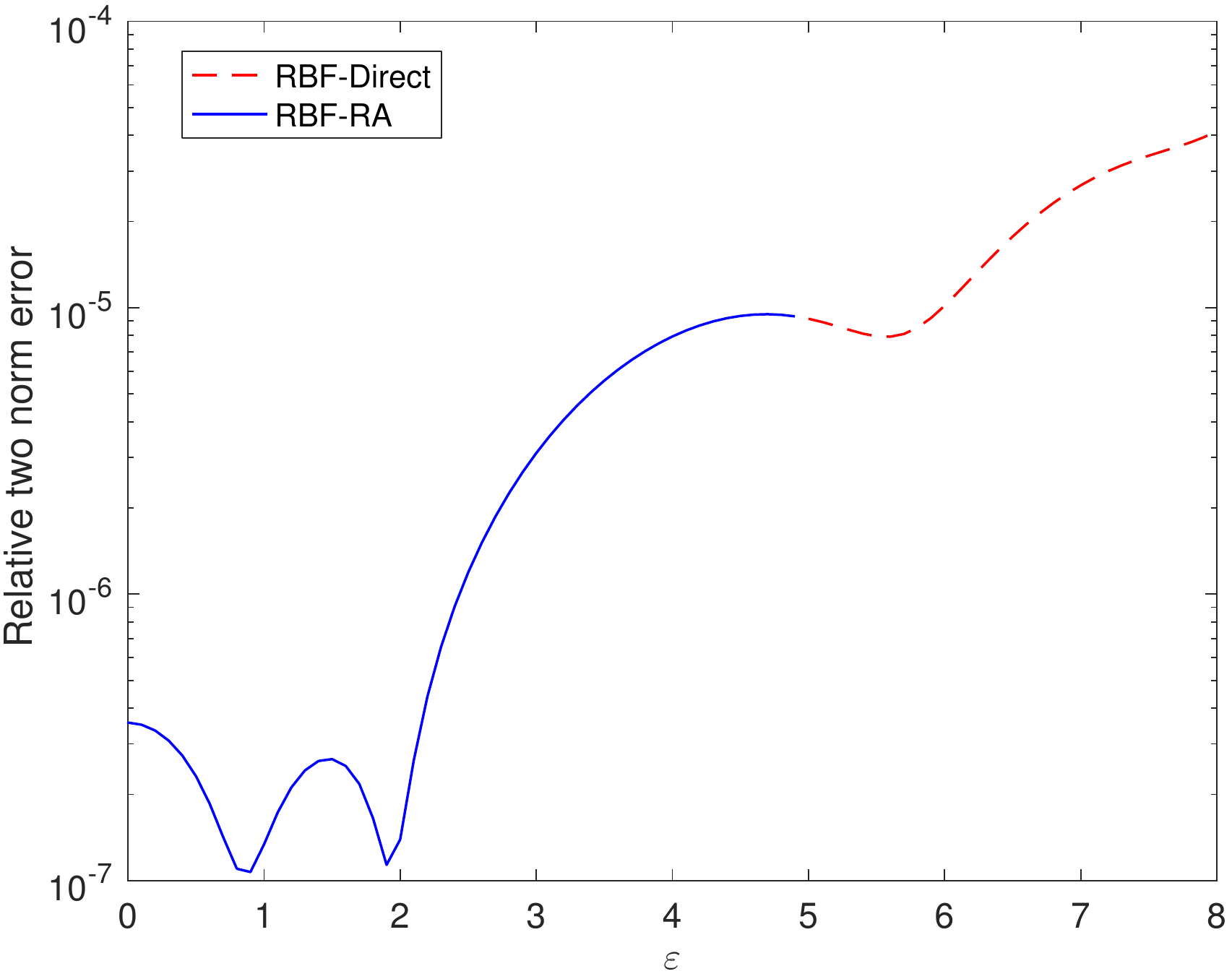} \\
(a) RBF-HFD weights & (b)  Solution to Poisson's equation
\end{tabular}
\caption{(a) Comparison of the errors in the computation of the RBF-HFD weights for one of the stencils from the node set shown in Figure \ref{fig:rbf_fd_illustration}.  The errors are measured against a multiprecision computation of the weights using 200 digits.  (b) Relative two-norm of the error in solving Poisson's equation \eqref{eq:Poisson} using RBF-HFD formulas as a function of $\ep$.  The dashed line marks the values of $\ep$ where RBF-Direct can be safely used to compute the RBF-HFD weights, while the solid line marks the values where RBF-RA is required to get an accurate result. The results in both plots were obtained using $N=45$ and $L=20$.\label{fig:rbf_hfd}}
\end{figure}

We first demonstrate the accuracy of the RBF-RA procedure for computing the RBF-HFD weights in the $\epz$ limit by comparing the computed weights for one of the stencils to the `exact' value of the weights computed using multiprecision arithmetic with $D=200$ digits.  The results are displayed in Figure \ref{fig:rbf_hfd} (a), together with the error in the weights computed with RBF-Direct in double precision arithmetic.  We can see from the figure that the RBF-RA method can compute the weights in a stable manner for the full range of $\ep$, with no loss in accuracy from the computation using RBF-Direct in the numerically `safe' region.  

Next, we use the RBF-HFD weights to numerically solve the Poisson equation \eqref{eq:Poisson} for various values of $\ep$ in the interval $[0,8]$.  
Letting $\underline{u}$ and $\underline{g}$ denote samples at the nodes $\Xi^{\text{int}}$ of the unknown solution and right hand side of \eqref{eq:Poisson}, respectively, the discretized version of \eqref{eq:Poisson} can be written as
\begin{equation*}
D\underline{u} = \dot{D}\underline{g}.
\end{equation*}
In this equation, $D$ and $\dot{D}$ are two sparse $P$-by-$P$ matrices with the $k^{\text{th}}$ row of $D$ containing the explicit weights $w_i$ in \eqref{eq:rbf_hfd} and the $k^{\text{th}}$ row of $\dot{D}$ containing the implicit weights $\tw_j$ in \eqref{eq:rbf_hfd}, for the $k^{\text{th}}$ node of $\Xi^{\text{int}}$.   To solve this system, we used BiCGSTAB with a zero-fill ILU-preconditioner and a tolerance on the relative residual of $10^{-10}$.   Figure \ref{fig:rbf_hfd} (b) shows the relative two-norm of the errors in the approximate solutions as a function of $\ep$.  Marked on this plot is the region where RBF-Direct becomes unstable, and RBF-RA is used.  We see that a reduction of the error by nearly two orders of magnitude is possible for $\ep$ values that are untouchable with RBF-Direct in double precision arithmetic.   We note that for all values of $\ep$, the solver required an average of 27.5 iterations and took 13.95 seconds of wall clock time using \matlab 2014b on a Linux workstation with 96 GB of RAM and dual 3.1GHz 8-core Intel Xeon processors.

\section{Concluding remarks}\label{sec:conclusion}
The present numerical tests demonstrate that the RBF-RA algorithm can be used effectively for stably computing RBF interpolants and RBF-FD/HFD formulas in the $\epz$ limit.  The method is more accurate and computationally efficient than the Contour-Pad\'e method, and more flexible than the RBF-QR and RBF-GA algorithms, in that operations such as differentiation can readily be applied directly to the RBFs instead of to the more complex bases that arise in these other two methods.  Its main disadvantages compared to these methods are (i) that it is not quite as efficient, and (ii) that it is more limited in the node sizes it can handle.  However, with the main target application of the algorithm being to compute RBF-FD and RBF-HFD formulas, these disadvantages are not serious concerns. The RBF-RA algorithm is also applicable to more general rational function approximation problems that involve vector-valued analytic functions with components that have the same singular points.  

A known issue with the RBF-CP, RBF-GA, and RBF-QR method is that the accuracy degrades if the nodes lie on a portion of a lower dimensional algebraic curve or surface (e.g.\ a cap on the unit sphere in $\mathbb{R}^3$).  This issue also extends to the RBF-RA method and will be a future topic of research.

\section*{Acknowledgments}

The work of Grady Wright was supported by National Science Foundation grants DMS-1160379 and ACI-1440638.

\appendix

\section{\matlab Code and examples}\label{appdx:code} This appendix contains first the script for two examples, followed by the functions \texttt{rbfhfd}, \texttt{vvra}, and \texttt{polyval2}. The function \texttt{vvra} (standing for \emph{vector valued rational approximation}) implements the general algorithm described in Section \ref{sec:vvra}.  The first example shows how to use the \texttt{vvra} function for RBF interpolation.  The output of this example is a plot showing the maximum difference between the interpolant and the target function \eqref{eq:TestFunction} as function of epsilon.  It also displays the minimum of these maximum differences and the corresponding value of $\ep$.  These numbers are $2.82\cdot 10^{-7}$ and $0.31$, respectively.  The second example shows how to use the \texttt{vvra} function for computing RBF-HFD weights.  The code uses the standard 19 node compact stencil for the 3D Laplacian~\cite{spotz_carey_1996} and shows that the RBF-HFD weights in the $\epz$ limit are the same (up to rounding errors) as the standard polynomials based weights.  The output of this example is the relative two norm difference between and flat limit RBF-HFD weights and the standard weights.  This number is $4.38\cdot 10^{-13}$.
\medskip

\begin{alltt}
\textcolor{comment}{%\% Example 1: Interpolation problem using GA kernel}
phi = @(e,r) exp(-(e*r).\^{}2); \textcolor{comment}{% Gaussian}
f = @(x,y) (1-(x.\^{}2+y.\^{}2)).*(sin(pi/2*(y-0.07))-0.5*cos(pi/2*(x+0.1)));
N = 60; hp = haltonset(2); y = 2*net(hp,N)-1;   \textcolor{comment}{% Nodes/centers}
M = 2*N; x = 3/2*net(scramble(hp,\textcolor{string}{'rr2'}),M)-3/4; \textcolor{comment}{% Eval points}
epsilon = linspace(0,1,101);
\textcolor{comment}{% Compute distances nodes/centers and eval points/centers}
D = @(x,y)hypot(bsxfun(@minus,x(:,1),y(:,1)'),bsxfun(@minus,x(:,2),y(:,2)'));
ryy = D(y,y); rxy = D(x,y);
\textcolor{comment}{% Determine the radius}
rad = fminbnd(@(e)norm(inv(phi(e,ryy)),inf)*norm(phi(1i*e,ryy),inf),0.1,20);
ryy = ryy*rad; \textcolor{comment}{% Re-scale the distances by the radius of the evaluation }
rxy = rxy*rad; \textcolor{comment}{% contour so that a unit radius can be used.}
\textcolor{comment}{% Compute the interpolant}
rbfinterp=@(ep) phi(ep,rxy)*(phi(ep,ryy)\(\backslash\)f(y(:,1),y(:,2)));
s = vvra(rbfinterp,epsilon/rad,1,64,64/4);  
\textcolor{comment}{% Compute the difference between s and f and plot the results}
error = max(abs(bsxfun(@minus,s,f(x(:,1),x(:,2)))));
semilogy(epsilon,error,\textcolor{string}{'x-'})
[minerr,pos] = min(error);
fprintf(\textcolor{string}{'Minimum error: \%1.2e, Epsion: \%1.2f\(\backslash\)n'},minerr,epsilon(pos));
\end{alltt}

\medskip

\begin{alltt}
\textcolor{comment}{%\% Example 2: RBF-HFD weights for the 3-D Laplacian using IQ kernel}
\textcolor{comment}{% Example for the standard 19 node stencil (6 implicit nodes) on lattice}
\textcolor{comment}{% Explicit nodes}
xhat = [[0,0,0];[-1,0,0];[1,0,0];[0,-1,0];[0,1,0];[0,0,-1];[0,0,1];\textcolor{keyword}{\underline{...}}
        [0,-1,-1];[0,-1,1];[0,1,-1];[0,1,1];[-1,0,-1];[-1,0,1];\textcolor{keyword}{\underline{...}}
        [1,0,-1];[1,0,1];[-1,-1,0];[-1,1,0];[1,-1,0];[1,1,0]];
N = size(xhat,1);
\textcolor{comment}{% Implicit nodes}
yhat = [[-1,0,0];[1,0,0];[0,-1,0];[0,1,0];[0,0,-1];[0,0,1]];
plot3(xhat(:,1),xhat(:,2),xhat(:,3),\textcolor{string}{'o'},yhat(:,1),yhat(:,2),yhat(:,3),\textcolor{string}{'.'})
\textcolor{comment}{% The IQ kernel and it's 3D Laplacian and 3D bi-harmonic}
phi = @(e,r) 1./(1 + (e*r).\^{}2);
dphi = @(e,r) 2*e\^{}2*(-3 + (e*r).\^{}2)./(1 + (e*r).\^{}2).\^{}3;
d2phi = @(e,r) 24*e\^{}4*(5 + (e*r).\^{}2.*(-10 + (e*r).\^{}2))./(1 + (e*r).\^{}2).\^{}5;
\textcolor{comment}{% Compute the radius}
rad = fzero(@(e)log10(cond(rbfhfd(e,xhat,yhat,phi,dphi,d2phi,1)))-6,[0.05,1]);
r = sqrt(D(xhat,xhat).\^{}2 + bsxfun(@minus,xhat(:,3),xhat(:,3)').\^{}2);
rad = min(rad,0.95/max(r(:)));
\textcolor{comment}{% Re-scale the stencil nodes so that a unit evaluation radius can be used.}
xhat = rad*xhat; yhat = rad*yhat;
\textcolor{comment}{% Compute the weights at various epsilon}
epsilon = linspace(0,rad,11);
w = vvra(@rbfhfd,epsilon/rad,1,64,64/4,xhat,yhat,phi,dphi,d2phi);  
\textcolor{comment}{% Undo the effects of re-scaling from the weights}
w(1:size(xhat,1),:) = w(1:size(xhat,1),:)*rad\^{}2;
\textcolor{comment}{% Compare the flat limit weights (epsilon=0) to the standard weights}
ws = [-8,2/3*ones(1,6),ones(1,12)/3,-ones(1,6)/6]';  \textcolor{comment}{% standard weights}
fprintf(\textcolor{string}{'Relative two norm difference: \%1.2e\(\backslash\)n'},norm(w(:,1)-ws)/norm(ws));
\end{alltt}

\medskip

\begin{alltt}
\textcolor{keyword}{function} w = rbfhfd(ep,x,xh,phi,dphi,d2phi,flag)
\textcolor{comment}{%RBFHFD  Computes the RBF-HFD weights for the 3D Laplacian.}
\textcolor{comment}{%}
\textcolor{comment}{%   w = rbfhfd(Epsilon,X,Xh,Phi,DPhi,D2Phi) computes the RBF-HFD weights }
\textcolor{comment}{%   at the given value of Epsilon and at the explicit stencil nodes X and}
\textcolor{comment}{%   implicit (Hermite) stencil nodes Xh.  Phi is a function handle for}
\textcolor{comment}{%   computing the kernel phi(ep,r) used for generating the weights (e.g. }
\textcolor{comment}{%   the Gaussian or inverse quadratic).  Dphi and D2phi are function}
\textcolor{comment}{%   handles for computing the 3D Laplacian and bi-harmonic of Phi,}
\textcolor{comment}{%   respectively.}
\textcolor{comment}{%}
\textcolor{comment}{%   A = rbfhfd(Epsilon,X,Xh,Phi,DPhi,D2Phi,flag) returns the matrix for }
\textcolor{comment}{%   computing the weights if the flag is non-zero, otherwise it returns}
\textcolor{comment}{%   just the weights as described above.}
\textcolor{keyword}{if} nargin == 6
    flag = 0;  \textcolor{comment}{% Return on the weights}
\textcolor{keyword}{end}
N = size(x,1);
x = [x;xh];
ov = ones(1,size(x,1));
\textcolor{comment}{% Compute the pairwise distances }
    r = sqrt((x(:,1)*ov-(x(:,1)*ov).').\^{}2 + (x(:,2)*ov-(x(:,2)*ov).').\^{}2 + \textcolor{keyword}{\underline{...}}
             (x(:,3)*ov-(x(:,3)*ov).').\^{}2);
temp = dphi(ep,r(1:N,N+1:end));  \textcolor{comment}{% Construct the weight matrix}
A = [[phi(ep,r(1:N,1:N)) temp];[temp.' d2phi(ep,r(N+1:\textcolor{keyword}{end},N+1:end))]];
\textcolor{comment}{% Determine what needs to be returned.}
\textcolor{keyword}{if} flag == 0
    w = A\(\backslash\)[dphi(ep,r(1:N,1));d2phi(ep,r(N+1:\textcolor{keyword}{end},1))];
\textcolor{keyword}{else}
    w = A;
\textcolor{keyword}{end}
end \textcolor{comment}{% end rbfhfd}
\end{alltt}

\medskip

\begin{alltt}
\textcolor{keyword}{function} [R,b] = vvra(myfun,epsilon,rad,K,n,varargin)
\textcolor{comment}{%VVRA   Vector-valued rational approximation (VVRA) valid near the}
\textcolor{comment}{%       origin to a vector-valued analytic function.}
\textcolor{comment}{%}
\textcolor{comment}{%   [R,b] = vvra(myfun,Epsilon,Rad,K,n) generates a vector-valued rational}
\textcolor{comment}{%   approximation to the vector-valued analytic function represented by}
\textcolor{comment}{%   myfun and evaluates it at the values in Epsilon. The inputs are}
\textcolor{comment}{%   described as follows:}
\textcolor{comment}{%   }
\textcolor{comment}{%      myfun is a function handle that, given a value of epsilon returns a}
\textcolor{comment}{%      column vector corresponding to each component of the function}
\textcolor{comment}{%      evaluated at epsilon.}
\textcolor{comment}{%}
\textcolor{comment}{%      Epsilon is an array of values such that |Epsilon|<=Rad where the}
\textcolor{comment}{%      rational approximation is to be evaluated.}
\textcolor{comment}{%}
\textcolor{comment}{%      Rad is a scalar representing the radius of the circle centered at}
\textcolor{comment}{%      the origin where it is numerically safe to evaluate myfun.}
\textcolor{comment}{%}
\textcolor{comment}{%      K is the number of points to evaluate myfun at on the contour for}
\textcolor{comment}{%      constructing the rational approximation. This number should be even.}
\textcolor{comment}{%}
\textcolor{comment}{%      2n is the degree of the common denominator to use in the }
\textcolor{comment}{%      approximation.}
\textcolor{comment}{% }
\textcolor{comment}{%   The columns in the output R represent the approximation to the}
\textcolor{comment}{%   components of myfun at each value in the array Epsilon. The optional}
\textcolor{comment}{%   output b contains the coefficients of the common denominator of the}
\textcolor{comment}{%   vector-valued rational approximation.}
\textcolor{comment}{%   }
\textcolor{comment}{%   [R,b] = vvra(myfun,Epsilon,Rad,K,n,T1,T2,...) is the same as above, but}
\textcolor{comment}{%   passes the optional arguments T1, T2, etc. to myfun, i.e.,}
\textcolor{comment}{%   feval(myfun,epsilon,T1,T2,...).}
K = K+mod(K,2);  \textcolor{comment}{% Force K to be even}
ang = pi/2*linspace(0,1,K+1)'; ang = ang(2:2:K);
ei = rad*exp(1i*ang);  \textcolor{comment}{% The evaluation points (all in first quadrant)}
m = K-n;               \textcolor{comment}{% Fix the degree of the numerator based on K and n}
W = feval(myfun,ei(1),varargin\{:\});
M = numel(W);
fv = zeros(K/2,1,M);
fv(1,1,:) = W;
\textcolor{keyword}{for} k = 2:K/2          \textcolor{comment}{% Loop over the evaluation points for F}
    W = feval(myfun,ei(k),varargin\{:\});
    fv(k,1,:) = W;
\textcolor{keyword}{end}
fmax = max(abs(fv),[],3);   \textcolor{comment}{% Find largest magnitude component for each k}
e = ei.\^{}2;                  \textcolor{comment}{% Calculate the E matrix}
E = e(:,ones(1,m)); E(:,1) = 1./fmax; E = cumprod(E,2);  \textcolor{comment}{% Scaled E-matrix}
f = E(:,1:n+1);             \textcolor{comment}{% Create the F-matrices and RHS}
F = f(:,:,ones(1,M)).*fv(:,ones(1,n+1),:);
g = F(:,1,:);
F = -F(:,2:n+1,:);
ER = [real(E);imag(E)]; \textcolor{comment}{% Separate E and F and the RHS g into real parts on }
FR = [real(F);imag(F)]; \textcolor{comment}{% top, and then the imag parts below}
gr = [real(g);imag(g)];
[Q,R] = qr(ER); QT = Q';  \textcolor{comment}{% Factorize ER into Q*R}
R = R(1:m,:);             \textcolor{comment}{% Remove bottom block of zeros from R}
\textcolor{keyword}{for} k = 1:M               \textcolor{comment}{% Update all FR matrices}
    FR(:,:,k) = QT*FR(:,:,k);
    gr(:,1,k) = QT*gr(:,1,k);
\textcolor{keyword}{end}
FT = FR(1:m,:,:);      \textcolor{comment}{% Separate F-blocks and g-blocks to the systems for}
FB = FR(m+1:K,:,:);    \textcolor{comment}{% determining the numerator and denominator coeffs.}
gt = gr(1:m,:,:);    
gb = gr(m+1:K,:,:);    
\textcolor{comment}{% Reshape these to be 2-D matrices}
FT = permute(FT,[1,3,2]); FT = reshape(FT,M*m,n);
FB = permute(FB,[1,3,2]); FB = reshape(FB,M*(K-m),n);
gt = permute(gt,[1,3,2]); gt = reshape(gt,M*m,1);
gb = permute(gb,[1,3,2]); gb = reshape(gb,M*(K-m),1);
b = FB\(\backslash\)gb;       \textcolor{comment}{% Obtain the coefficients of the denominator}
v = gt-FT*b; V = reshape(v,m,M);
a = (R\(\backslash\)V);       \textcolor{comment}{% Obtain the coefficients of the numerators}
\textcolor{comment}{% Evaluate the rational approximations}
R = zeros(M,length(epsilon));
b = [1;b];
denomval = polyval2(b,epsilon);
\textcolor{keyword}{for} ii = 1:M
  R(ii,:) = (polyval2(a(:,ii),epsilon)./denomval);
\textcolor{keyword}{end}
end  \textcolor{comment}{% End vvra}
\end{alltt}

\medskip

\begin{alltt}
\textcolor{keyword}{function} y = polyval2(p,x)
\textcolor{comment}{%POLYVAL2 Evaluates the even polynomial}
\textcolor{comment}{%         Y = P(1) + P(2)*X\^{}2 + ... + P(N)*X\^{}(2(N-1)) + P(N+1)*X\^{}2N}
\textcolor{comment}{% If X is a matrix or vector, the polynomial is evaluated at all}
\textcolor{comment}{% points in X (this is unlike the polyval function of matlab)}
y = zeros(size(x)); x = x.\^{}2;
\textcolor{keyword}{for} j=length(p):-1:1
    y = x.*y + p(j);
\textcolor{keyword}{end}
end  \textcolor{comment}{% End polyval2}
\end{alltt}

\section{Choosing the radius $\ep_R$ of the evaluation contour}\label{appdx:contours}
We propose two different strategies for choosing the radius $\ep_R$ depending on the type of kernel being used in the application.  As discussed in Section \ref{sec:interpolation}, entire positive definite kernels, like GA, grow without bound as one moves out in the complex plane away from the real $\ep$-axis.  This growth limits the radius of the contour.  The strategy we propose for entire kernels is to choose $\ep_R$ based on an estimate of where the ill-conditioning from the interpolation problem is overtaken by ill-conditioning from the growth of the kernel on the imaginary $\ep$-axis.  Such an estimate is difficult to make from the condition number of the interpolation or RBF-FD/HFD weight matrix $A(\ep)$ because of singularities that can occur in $A(\ep)^{-1}$ for imaginary $\ep$ (cf. Figure \ref{fig:cond_surface}).  Instead, we have found that a good measure for the ill-conditioning that is not effected by singularities when $\ep=i\beta$ ($\beta \in \mathbb{R}$) is
\begin{align*}
\widetilde{\sigma}_{\infty}(A(\beta)) = \|A(i\beta)\|_{\infty} \|A(\beta)^{-1}\|_{\infty}.
\end{align*}
The first factor on the right of this expression captures the growth from the kernel, while the second factor captures the ill-conditioning from the interpolation problem (since this is just the standard interpolation matrix for real shape parameters).  To find where the transition between the types of ill-conditioning occurs, we thus set $\ep_R$ equal to an approximate minimum of $\log\widetilde{\sigma}_{\infty}(A(\beta))$.   This is the approach used in the interpolation example in the first appendix.

As discussed in Section \ref{sec:interpolation}, kernels that have singular points, such as IQ, IMQ, and MQ, will lead to a multitude of singular points in the RBF vector-valued rational function, as each entry of $A(\ep)$ will have two singular points.   The evaluation contour that is used should altogether avoid these singular points.  The closest singular points that comes directly from the kernels is located on the positive and negative imaginary axis at the inverse of the largest distance between the nodes.  The radius $\ep_R$ thus needs to be chosen smaller than this value.  In our applications we have used
\begin{align}
\ep_R = 0.95 \left(\max_{1\leq i,j, \leq N} \| \vxc_i - \vxc_j \|\right)^{-1}, \label{eq:epR2}
\end{align}
where $\{\vxc_{i}\}_{i=1}^N$ are the interpolation (or RBF-FD stencil) nodes.  For smaller problems, we have found that this sometimes results in a radius that is larger than necessary to achieve accurate results.  In these case we often find it is better to use a smaller radius, which in turn allows the number of evaluation points $K$ to be smaller without diminishing the accuracy.  To address this issue, we choose $\ep_R$ to be the minimum between \eqref{eq:epR2} and the approximate real value of $\ep$ where the condition number of $A(\ep)$ is equal to $10^6$.  This is the approach used in the RBF-HFD example in the first appendix.

\section*{References}
\bibliography{rbf_ra_refs}

\end{document}